\documentclass[a4paper,reqno,11pt]{amsart}

\usepackage[utf8]{inputenc}
\usepackage[T1]{fontenc}
\usepackage[english]{babel}
\usepackage{lmodern}
\usepackage[margin=1.0in]{geometry}
\usepackage{microtype}

\usepackage{graphicx}
\usepackage{amsmath,amssymb,amsthm}
\usepackage{mathtools}
\usepackage{bm}
\usepackage[dvipsnames]{xcolor}
\usepackage{enumitem}
\usepackage{booktabs}

\usepackage[hidelinks]{hyperref}

\allowdisplaybreaks
\numberwithin{equation}{section}

\theoremstyle{plain}
\newtheorem{theorem}{Theorem}

\newtheorem{lemma}{Lemma}
\newtheorem{corollary}{Corollary}

\newtheorem{proposition}{Proposition}

\theoremstyle{definition}

\newtheorem{remark}[theorem]{Remark}

\theoremstyle{remark}

\DeclareMathOperator{\dist}{dist}

\newcommand{\R}{{\mathbb{R}}}

\newcommand{\Sphere}{\mathbb S}
\newcommand{\dd}{\,\mathrm d}
\newcommand{\ind}{\mathbf 1}

\DeclareMathOperator{\Int}{int}

\DeclareMathOperator*{\argmax}{arg\,max}

\definecolor{darkred}{RGB}{139,0,0}

\title{Ehrhard-Type Inequalities for Radial Measures}
\author{Tomasz Byczkowski}

\address{Department of Mathematics,
Faculty of Pure and Applied Mathematics,
Wroc{\l}aw University of Science and Technology,
ul. Wybrze{\.z}e Wyspia{\'n}skiego 27,
50-370 Wroc{\l}aw, Poland}

\email{tomasz.byczkowski@gmail.com}

\author{Jacek Ma{\l}ecki}

\thanks{Corresponding author: Jacek Ma{\l}ecki.
E-mail:
\href{mailto:jacek.malecki@pwr.edu.pl}
{jacek.malecki@pwr.edu.pl}.}

\address{Department of Mathematics,
Faculty of Pure and Applied Mathematics,
Wroc{\l}aw University of Science and Technology,
ul. Wybrze{\.z}e Wyspia{\'n}skiego 27,
50-370 Wroc{\l}aw, Poland}

\email{jacek.malecki@pwr.edu.pl}

\author{Tomasz Żak}

\address{Department of Mathematics,
Faculty of Pure and Applied Mathematics,
Wroc{\l}aw University of Science and Technology,
ul. Wybrze{\.z}e Wyspia{\'n}skiego 27,
50-370 Wroc{\l}aw, Poland}

\email{tomasz.zak@pwr.edu.pl}

\date{\today}

\subjclass[2020]{60E15, 52A40, 60E07}

\keywords{Ehrhard inequality; radial probability measures; sectional symmetrization; reciprocal logarithmic slope; Student mixtures; stable distributions.}

\begin{document}

\begin{abstract}
We study Ehrhard-type inequalities for convex sets under rotationally invariant probability measures on $\R^n$, $n\ge2$, described by a radial density profile $\varphi$. We give sufficient conditions on the profile that apply to Gaussian and multivariate Student measures, including the Cauchy measure, and to a family of mixtures of two Student measures. We also prove necessary conditions on the profile and its tails. Together, these conditions determine the exact parameter range for this mixture family and characterize the rotationally invariant stable measures satisfying the inequality.
\end{abstract}

\maketitle

\section{Introduction and main results}\label{sec:intro}
In 1983, Ehrhard \cite{Ehrhard} proved a Gaussian counterpart of the Brunn--Minkowski inequality for Lebesgue measure. More precisely, the map $K\mapsto\Phi^{-1}(\gamma_n(K))$ is concave under Minkowski combinations of convex Borel sets in $\R^n$, where $\gamma_n$ is the standard Gaussian measure and $\Phi$ is the distribution function of $N(0,1)$. His proof used Gaussian $k$-symmetrizations, $k=1,\ldots,n$, and their concavity properties. Lata{\l}a \cite{Latala} extended the inequality to the case where only one of the two Borel sets is convex. Borell \cite{Borell} later proved it for arbitrary Borel sets using the heat flow and the Feynman-Kac formula.

\smallskip

The analogous question for the isotropic Cauchy measure was studied by Byczkowski and \.{Z}ak \cite{BZ}. They introduced a Cauchy $1$-symmetrization in $\R^n$ and proved joint concavity of the corresponding section map. Byczkowski, Ma{\l}ecki and \.{Z}ak \cite{BMZ} then proved an analogue of Ehrhard's inequality for axis-parallel planar rectangles symmetric about a common coordinate axis. They also conjectured that the inequality holds for all convex Borel sets in $\R^n$. Huang, Tang and Wang \cite{HTW} subsequently proved this conjecture for the isotropic Cauchy measure and established the corresponding inequality for a family of Cauchy-type measures in every dimension. Their proof combines a variational argument with joint concavity of the section rearrangement.

\smallskip

In this paper, we give sufficient conditions on the radial density profile for an Ehrhard-type inequality in $\R^n$. These conditions apply to Gaussian and multivariate Student measures, including the Cauchy measure, and to a family of mixtures of two Student measures. We also derive necessary conditions on the profile and its tails. Combining the sufficient and necessary conditions, we determine the exact parameter range for this mixture family and characterize the rotationally invariant stable measures satisfying the inequality.

\smallskip

We first give sufficient conditions directly on the radial density profile. They apply also to measures outside the Gaussian and Student families. We then derive necessary conditions directly from the Ehrhard-type inequality, without imposing the sufficient assumptions. These necessary conditions exclude gaps in the radial support and bounded support. They also restrict the possible power-law tails and give exact parameter ranges for the new mixture family and for isotropic stable measures.

Throughout, $n\ge2$ and
    \begin{equation}
        \label{eq:R}
        \begin{gathered}
            Z_\varphi=\int_{\R^n}\varphi(|x|^2)\dd x,\qquad 0<Z_\varphi<\infty,\\
            \rho(x)=Z_\varphi^{-1}\varphi(|x|^2),\qquad \dd\mu_\varphi(x)=\rho(x)\dd x.
        \end{gathered}
    \end{equation}
Thus the argument of $\varphi$ is the square of the distance from the origin. Unless stated otherwise, $\varphi$ is continuous and nonnegative. When $\varphi$ is strictly positive, we write $F_\varphi$ for the distribution function of the first coordinate and $\Psi_\varphi$ for its inverse:
    \begin{equation*}
        F_\varphi(t)=\mu_\varphi\{x:x_1\le t\}, \qquad \Psi_\varphi=F_\varphi^{-1}:(0,1)\to\mathbb R.
    \end{equation*}
By rotational invariance, $F_\varphi$ is also the distribution function of $u\cdot X$ for every unit vector $u$, where $X$ has law $\mu_\varphi$. It is not the distribution function of $|X|$, and it is generally different from the distribution function on a normalized one-dimensional section. When $\varphi$ is fixed, we also write $\mu=\mu_\varphi$, $F=F_\varphi$, and $\Psi=F_\varphi^{-1}$. The property studied in this paper is the following inequality:
    \begin{equation}
        \label{eq:E}\tag{E}
        \Psi_\varphi\!\left(\mu_\varphi((1-\lambda)K+\lambda L)\right) \ge (1-\lambda)\Psi_\varphi(\mu_\varphi(K)) +\lambda\Psi_\varphi(\mu_\varphi(L)),
    \end{equation}
for convex Borel sets $K,L$ with measures in $(0,1)$ and $\lambda\in(0,1)$. We define $F_\varphi^{-1}(0)$ and $F_\varphi^{-1}(1)$ by the corresponding one-sided limits. For a strictly positive profile, the support is $\R^n$, so these values are $-\infty$ and $+\infty$, respectively. In particular, the left-hand side of \eqref{eq:E} is $+\infty$ when the Minkowski combination has measure $1$.

\subsection{A sufficient criterion on the squared-radius profile}

If $\varphi$ is differentiable, $\varphi(s)>0$ and $\varphi'(s)<0$ for every $s\ge0$, let
    \begin{equation}
        \label{eq:w}
        w(s)=-\frac{\varphi(s)}{2\varphi'(s)},\qquad s\ge0.
    \end{equation}
Our first theorem is expressed solely through this scalar function.

\begin{theorem}[Reciprocal-slope criterion]
    \label{thm:main}
    Let $n\ge2$ and $\varphi\in C^3([0,\infty))$ be strictly positive, with $\varphi'(s)<0$ for all $s\ge0$ and
        \begin{equation*}
            \int_0^\infty r^{n-1}\varphi(r^2)\dd r<\infty.
        \end{equation*}
    If $w$ from \eqref{eq:w} satisfies 
        \begin{equation*}
            w'(s)\ge0,\qquad w''(s)\ge0,\qquad w(s)-s w'(s)\ge0\quad(s\ge0),
        \end{equation*}
    then $\mu_\varphi$ satisfies \eqref{eq:E} for all convex Borel sets $K,L\subseteq\mathbb R^n$ with $0<\mu_\varphi(K),\mu_\varphi(L)<1$.
\end{theorem}

The criterion applies to the Gaussian profile $\varphi(s)=e^{-\beta s}$, for which $w$ is constant, and to the Student profiles $\varphi(s)=(\tau+s)^{-p}$, for which $w$ is affine. It also applies when $w$ is not affine, for example to
    \begin{equation}
        \label{eq:preview-mixture}
        \varphi_\varepsilon(s) = (\tau+s)^{-p}+\varepsilon(\tau+s)^{-p-1}, \qquad p>\frac n2,\quad\tau>0, \quad 0\le\varepsilon\le p\tau.
    \end{equation}
For $\varepsilon>0$ this profile does not give a Gaussian measure or a single rescaled Student measure. Section~\ref{sec:examples} verifies the sufficient criterion for this family. Corollary~\ref{cor:necessary-mixture-iff} below shows that the parameter range in \eqref{eq:preview-mixture} is exact.

\subsection{Necessary profile and support conditions}

We now ask what the inequality \eqref{eq:E} implies about the profile. We do not assume the conditions on $w$ in this part. In the next theorem we also allow $\varphi$ to vanish. For such profiles we use the generalized inverse
        \begin{equation*}
            F_\varphi^{-1}(p)=\inf\{t\in\R:F_\varphi(t)\ge p\}, \qquad 0<p<1,
        \end{equation*}
and define its values at $0$ and $1$ by one-sided limits.

\begin{theorem}[Necessary positivity, monotonicity, and full support]
    \label{thm:necessary-full-support}
    Let $n\ge2$, and let $\varphi:[0,\infty)\to[0,\infty)$ be continuous, not identically zero, with $0<Z_\varphi<\infty$. If \eqref{eq:E} holds for all convex sets $K,L$ with $0<\mu_\varphi(K),\mu_\varphi(L)<1$, then
        \begin{equation}
            \label{eq:necessary-support}
            \operatorname{supp}\mu_\varphi=\mathbb R^n,\qquad \varphi(s)>0\ (s\ge0),\qquad \varphi\text{ is nonincreasing}.
        \end{equation}
\end{theorem}

The proof uses a test with two small balls and an infinitesimal perturbation of half-spaces. In particular, it shows that for $Z=(X,Y)\sim\mu_\varphi$, with $X\in\mathbb R$ and $Y\in\mathbb R^{n-1}$, the conditional mean
    \begin{equation*}
        t\longmapsto\mathbb E[|Y|\mid X=t]
    \end{equation*}
must be even and convex on the interior of the marginal support. This property is incompatible with any bounded radial support. Thus \eqref{eq:E} requires $\varphi$ to be strictly positive and nonincreasing. The assumption $\varphi'(s)<0$ in Theorem~\ref{thm:main} is stronger. We use it to define $w$; we do not claim that \eqref{eq:E} implies $\varphi'(s)<0$ at every point.

For measures with power-law tails, we obtain another necessary condition by applying \eqref{eq:E} to balls of the form $x+\varepsilon rB^n$ and letting $\varepsilon\downarrow0$. Here $B^n=\{z\in\R^n:|z|\le1\}$ is the Euclidean unit ball. %\tz{skoro ''unit ball'', to czemu $B_2$ a nie $B_1$?} 
The balls shrink to their centers, so their measures tend to zero. The tail assumption below describes $F_\varphi^{-1}(p)$ as $p\downarrow0$ and is used in this limiting argument.

\begin{theorem}[Necessary convexity condition from small balls]
    \label{prop:necessary-small-balls}
    Let $\varphi$ be continuous and strictly positive, with $0<Z_\varphi<\infty$. Suppose that its marginal satisfies
        \begin{equation}
            \label{eq:necessary-marginal-tail}
            1-F_\varphi(t)\sim A t^{-\eta} \qquad(t\to\infty),\qquad A,\eta>0.
        \end{equation}
    If \eqref{eq:E} holds, then
        \begin{equation}
            \label{eq:necessary-density-convexity}
            x\longmapsto\rho(x)^{-1/(n+\eta)} \text{ is convex on }\mathbb R^n,
        \end{equation}
    Here $\rho$ is the density defined in \eqref{eq:R}.
\end{theorem}

An immediately applicable consequence uses only the second-order asymptotics of the squared-radius profile.

\begin{theorem}[Second-order profile-tail obstruction]
    \label{thm:necessary-profile-tail}
    Let $\varphi$ be continuous, strictly positive and integrable as in \eqref{eq:R}, and assume that
        \begin{equation}
            \label{eq:necessary-profile-asymptotic}
            \varphi(s)=C s^{-(n+\eta)/2} \left(1+b s^{-\theta}+o(s^{-\theta})\right) \qquad(s\to\infty),
        \end{equation}
    where $C,\eta,\theta>0$ and $b\ne0$. A necessary condition for \eqref{eq:E} is
        \begin{equation*}
            b(2\theta-1)\le0.
        \end{equation*}
    Equivalently, if $b(2\theta-1)>0$, then \eqref{eq:E} fails for some convex sets.
\end{theorem}

When $2\theta=1$, this test gives no conclusion; further terms of the expansion may be needed. The proof does not differentiate the remainder and does not assume joint concavity of the function $g$ defined in \eqref{eq:threshold}.

\subsection{Sharp parameter classifications}

For the family \eqref{eq:preview-mixture}, the sufficient and necessary conditions give the same parameter range.

\begin{corollary}[Sharp range for adjacent-power mixtures]
    \label{cor:necessary-mixture-iff}
    For $n\ge2$, $p>n/2$, $\tau>0$, and $\varepsilon\ge0$, the probability measure with squared-radius profile
        \begin{equation*}
            \varphi_\varepsilon(s) = (\tau+s)^{-p}+\varepsilon(\tau+s)^{-p-1}
        \end{equation*}
    satisfies \eqref{eq:E} for all convex sets if and only if
        \begin{equation*}
            \quad0\le\varepsilon\le p\tau.
        \end{equation*}
\end{corollary}

Theorem~\ref{thm:main} proves sufficiency. For necessity, we use Theorem~\ref{thm:necessary-profile-tail} with $\eta=2p-n$, $\theta=1$, and $b=\varepsilon-p\tau$.

The same necessary condition determines which rotationally invariant stable measures satisfy \eqref{eq:E}. For $0<\alpha<2$, $\alpha\ne1$, the first coefficient in the density expansion is positive, while the second has the sign of $\alpha-1$. Written in the form \eqref{eq:necessary-profile-asymptotic}, this expansion has $\eta=\alpha$ and $\theta=\alpha/2$. Hence $b(2\theta-1)=b(\alpha-1)>0$, which excludes \eqref{eq:E}. Together with the Cauchy and Gaussian cases, this gives Corollary~\ref{thm:stable-classification} below.

\begin{corollary}[Classification of isotropic stable laws]
    \label{thm:stable-classification}
    Let $n\ge2$, $0<\alpha\le2$, and let $\mu_{n,\alpha}$ be the rotationally invariant stable probability measure with
        \begin{equation}
            \label{eq:stable-characteristic}
            \widehat\mu_{n,\alpha}(\xi)=e^{-|\xi|^\alpha}.
        \end{equation}
    Writing $F_\alpha$ for the distribution function of its first coordinate, \eqref{eq:E} with $\mu=\mu_{n,\alpha}$ and $\Psi=F_\alpha^{-1}$ holds for all convex sets if and only if
        \begin{equation}
            \label{eq:stable-iff}
            \alpha=1\ \text{or}\ \alpha=2.
        \end{equation}
    For every $\alpha\in(0,2)\setminus\{1\}$, the inequality fails at $\lambda=1/2$ for some pair of bounded boxes with sides parallel to the coordinate axes. Both boxes are symmetric with respect to every coordinate hyperplane.
\end{corollary}

The Gaussian and Cauchy cases follow from the sufficient criterion. The necessary tail condition excludes every other index. To obtain counterexamples among bounded boxes, we first prove failure for slabs using the one-dimensional marginal, and then truncate the slabs.

\subsection{Method and organization.}

We first prove that the endpoint of the half-line replacing a section $(a,b)$ is jointly concave in the position of the section and in $a,b$. The direction of the section is fixed in this step. We then use symmetrization to prove \eqref{eq:E} for convex sets. The proof uses the ideas of Ehrhard and Byczkowski-\.{Z}ak and is related to the argument in \cite{HTW}; we do not use the multidimensional Student inequality as an input. Section~\ref{sec:analytic} proves the required concavity from the conditions on $w$. Section~\ref{sec:geometry} completes the proof by symmetrization. Section~\ref{sec:examples} gives the examples. Section~\ref{sec:necessary} proves the necessary conditions and the exact range for the mixture family. Finally, Section~\ref{sec:stable-classification} proves the stable classification and constructs counterexamples among bounded boxes.

\section{Joint concavity of the section map}
\label{sec:analytic}

In this section, we replace an interval on each line by a left half-line with the same integral of the density, and study the endpoint of that half-line. We prove that this endpoint depends concavely on the position of the line and on the two interval endpoints; Section~\ref{sec:geometry} will use this property to show that symmetrization sends convex sets to convex sets. Recall the density $\rho(x)=Z_\varphi^{-1}\varphi(|x|^2)$ from \eqref{eq:R}. Fix a unit vector $v\in\Sphere^{n-1}$ and write $x=y+tv$, where $y\in v^\perp$ gives the position of the line parallel to $v$. Set
    \begin{equation*}
        H_y(s)=\int_{-\infty}^s\rho(y+tv)\dd t, \qquad M(y)=\int_{\mathbb R}\rho(y+tv)\dd t.
    \end{equation*}
For finite $a<b$, let $g(y,a,b)$ be the unique threshold satisfying
    \begin{equation}
        \label{eq:threshold}
        H_y(g(y,a,b))=H_y(b)-H_y(a).
    \end{equation}
The integral of $t\mapsto\rho(y+tv)$ is $M(y)$, which need not equal $1$. We denote the normalized density and distribution function by
    \begin{equation*}
        \widetilde f_y(t)=\frac{\rho(y+tv)}{M(y)}, \qquad \widetilde F_y(s)=\frac{H_y(s)}   {M(y)}.
    \end{equation*}
Here $y$ specifies the line, and integration is with respect to $t$. Since $M(y)>0$, equation \eqref{eq:threshold} is equivalent to $\widetilde F_y(g)=\widetilde F_y(b)-\widetilde F_y(a)$. Thus normalization does not change $g$. We will prove that the assumptions of Theorem~\ref{thm:main} imply
    \begin{equation}
        \label{eq:JC}\tag{JC}
        (y,a,b)\longmapsto g(y,a,b) \quad\text{is jointly concave on }v^\perp\times\{a<b\}.
    \end{equation}
    Equivalently, for $0<\lambda<1$, $y_0,y_1\in v^\perp$, and $a_0<b_0$, $a_1<b_1$, we have
    \begin{align}
        \label{eq:JC-expanded}
        &g((1-\lambda)y_0+\lambda y_1, (1-\lambda)a_0+\lambda a_1, (1-\lambda)b_0+\lambda b_1)\notag\\
        &\hspace{15mm}\ge(1-\lambda)g(y_0,a_0,b_0)+\lambda g(y_1,a_1,b_1).
    \end{align}
The direction $v$ is fixed in \eqref{eq:JC}; the variables are the position $y$ of the line and the two interval endpoints $a,b$. By radiality, the same assertion then holds for every unit direction $v$. We prove it using a comparison of $2\times2$ matrices.

\subsection{A matrix criterion for normalized one-dimensional densities}

In the following lemma, $y\in\R$ is a parameter. For every fixed $y$, the function $x\mapsto f(y,x)$ is a probability density on $\R$, and
    \begin{equation*}
        \int_{\R}f(y,x)\dd x=1,\qquad f(y,-x)=f(y,x).
    \end{equation*}
The functions $f(y,\cdot)$, indexed by $y\in\R$, form a family of one-dimensional probability densities. The two-variable function $f$ describes this family; it is not a joint probability density on $\R^2$.

The parameter $y$ is allowed to vary, and derivatives with respect to $y$ are taken after each section has been normalized. We assume that $f>0$, $f_x(y,x)<0$ for $x>0$, and $f_{xx}(y,0)<0$. We also assume that $f$ and the upper-tail function
    \begin{equation*}
        P(y,x)=\int_x^\infty f(y,s)\dd s
    \end{equation*}
are $C^2$ jointly in $(y,x)$. In the application to the radial profile, we will normalize the sections and check these assumptions. Symmetry in $x$ and normalization give
    \begin{equation*}
        P(y,0)=\int_0^\infty f(y,s)\dd s =\frac12\int_{\R}f(y,s)\dd s=\frac12 \qquad\text{for every }y.
    \end{equation*}
For $x>0$ set
    \begin{equation}
        \label{eq:rs-tail-matrix}
        \begin{split}
            J(y,x)&=-\frac{f(y,x)^2}{f_x(y,x)}>0, \qquad k(y,x)=P_y(y,x)+\frac{f(y,x)f_y(y,x)}{f_x(y,x)},\\
            \ell(y,x)&=-P_{yy}(y,x)-\frac{f_y(y,x)^2}{f_x(y,x)}, \qquad \mathcal M(y,x)=
                \begin{pmatrix}J(y,x)&k(y,x)\\k(y,x)&\ell(y,x)\end{pmatrix}.
        \end{split}
    \end{equation}
For fixed $y$, the function $x\mapsto P(y,x)$ is continuous and strictly decreasing on $[0,\infty)$, with $P(y,0)=1/2$ and $P(y,x)\to0$ as $x\to\infty$. Hence, for every $p\in(0,1/2)$, there is a unique $x_y(p)>0$ such that $P(y,x_y(p))=p$. We write
    \begin{equation*}
        \mathcal M_y(p)=\mathcal M(y,x_y(p)).
    \end{equation*}
For example, consider $\rho(y,x)=\pi^{-1}e^{-y^2-x^2}$, the density of the centered Gaussian law with covariance matrix $\tfrac12 I_2$. The integral over the line is $M(y)=\pi^{-1/2}e^{-y^2}$. After division by $M(y)$, we obtain
    \begin{equation*}
        f(y,x)=\frac{e^{-x^2}}{\sqrt\pi},\qquad P(y,x)=1-\Phi(\sqrt2\,x),\qquad x_y(p)=\frac{\Phi^{-1}(1-p)}{\sqrt2},
    \end{equation*}
where $\Phi$ is the standard Gaussian distribution function. In this example the normalized density and $x_y(p)$ do not depend on $y$. The unnormalized integral over $[0,\infty)$ is $M(y)/2$, whereas $P(y,0)=1/2$.

For real symmetric $2\times2$ matrices $A$ and $B$, we write $A\succeq B$ when $A-B$ is positive semidefinite, that is,
    \begin{equation*}
        \xi^T(A-B)\xi\ge0\qquad\text{for every }\xi\in\R^2.
    \end{equation*}
Thus a matrix-valued function is nondecreasing in this order when each of its quadratic forms is nondecreasing. In \eqref{eq:rs-matrix-monotonicity} below, the assumption means
    \begin{equation*}
        \frac{\mathcal M_y(q)}q-\frac{\mathcal M_y(p)}p\succeq0 \qquad(0<p<q<1/2).
    \end{equation*}
It does not mean that each matrix entry is nondecreasing.

\begin{lemma}[Matrix tail criterion]
    \label{lem:rs-matrix-tail}
    Under the preceding assumptions, suppose that, for every fixed $y$,
        \begin{equation}
            \label{eq:rs-matrix-monotonicity}
            p\longmapsto\frac{\mathcal M_y(p)}p \quad\text{is nondecreasing on }(0,1/2)
        \end{equation}
    in the positive-semidefinite order. Then the function $g$ defined by
        \begin{equation*}
            \int_{-\infty}^{g(y,a,b)}f(y,s)\,ds =\int_a^b f(y,s)\,ds,\qquad a<b,
        \end{equation*}
    is jointly concave in $(y,a,b)$.
\end{lemma}

\begin{proof}
We first record two consequences of the assumed matrix monotonicity. We then use them to show that the second derivative of $g$ along every line in its domain is nonpositive. Fix $y$ for the moment. Matrix monotonicity implies superadditivity:
    \begin{equation}
        \label{eq:rs-superadditivity}
        \mathcal M_y(p+q)\succeq \mathcal M_y(p)+\mathcal M_y(q), \qquad p,q>0,\quad p+q<1/2.
    \end{equation}
Indeed, compare each of $\mathcal M_y(p)/p$ and $\mathcal M_y(q)/q$ with $\mathcal M_y(p+q)/(p+q)$ and add.

We next find a vector $e_y\ne0$, depending on $y$ but not on $p$, such that the quadratic form of $\mathcal M_y(p)$ evaluated at $e_y$ is nonpositive for every $0<p<1/2$. Put
    \begin{equation*}
        \kappa(y)=-\frac{f_y(y,0)}{f(y,0)},\qquad e_y=(-\kappa(y),1)^T.
    \end{equation*}
A direct calculation from \eqref{eq:rs-tail-matrix} gives
    \begin{equation}
        \label{eq:rs-rank-one-form}
        \begin{aligned}
            e_y^T\mathcal M(y,x)e_y &= J(y,x) \left(\kappa(y)+\frac{f_y(y,x)}{f(y,x)}\right)^2\\ &\quad -2\kappa(y)P_y(y,x)-P_{yy}(y,x).
        \end{aligned}
    \end{equation}
We show that each term on the right tends to zero as $x\downarrow0$, with $y$ fixed.

Symmetry gives $f_x(y,0)=0$. Since $f$ is $C^2$ and $f_{xx}(y,0)<0$,
    \begin{equation*}
        f_x(y,x)=f_{xx}(y,0)x+o(x), \qquad J(y,x)= -\frac{f(y,0)^2}{f_{xx}(y,0)} \frac{1+o(1)}{x} = O(x^{-1}).
    \end{equation*}
Differentiating the identity $f_x(y,0)=0$ with respect to $y$ also gives $f_{yx}(y,0)=0$. Thus the derivative in $x$ of $f_y/f$ vanishes at $x=0$, and
    \begin{equation*}
        \kappa(y)+\frac{f_y(y,x)}{f(y,x)} = \frac{f_y(y,x)}{f(y,x)} -\frac{f_y(y,0)}{f(y,0)} =o(x).
    \end{equation*}
It follows that the first term in \eqref{eq:rs-rank-one-form} is $o(x)$. Moreover, $P(y,0)=1/2$ for every $y$, so $P_y(y,0)=P_{yy}(y,0)=0$. Continuity of these derivatives shows that the other two terms also tend to zero.

Finally, $x_y(p)\to0$ as $p\uparrow1/2$, while $p$ tends to the nonzero value $1/2$. Therefore
    \begin{equation*}
        \lim_{p\uparrow1/2} e_y^T\frac{\mathcal M_y(p)}p e_y=0.
    \end{equation*}
    
Together with \eqref{eq:rs-matrix-monotonicity}, this gives
    \begin{equation}
        \label{eq:rs-common-direction}
        e_y^T\mathcal M_y(p)e_y\le0, \qquad 0<p<1/2.
    \end{equation}

We will express derivatives of $g$ in terms of derivatives of the tail probabilities $P(y,x)$. The following quadratic form will be useful. For a symmetric matrix $M=\left(\begin{smallmatrix}J&k\\k&\ell\end{smallmatrix}\right)$ with $J>0$, define
    \begin{equation}
        \label{eq:rs-quadratic-form}
        \mathcal B_M(\eta,\sigma) =\frac{(\sigma-k\eta)^2}{J}-\ell\eta^2 =\sup_{q\in\mathbb R} \left\{2q\sigma-Jq^2-2kq\eta-\ell\eta^2\right\}.
    \end{equation}
Two consequences will be used:
    \begin{align}
        M\succeq N&\quad\Longrightarrow\quad \mathcal B_M(\eta,\sigma)\le\mathcal B_N(\eta,\sigma), \label{eq:rs-order-B}\\
        \mathcal B_{M_1+M_2}(\eta,\sigma_1+\sigma_2) &\le\mathcal B_{M_1}(\eta,\sigma_1) +\mathcal B_{M_2}(\eta,\sigma_2). \label{eq:rs-sum-B}
    \end{align}
In \eqref{eq:rs-order-B}, the $(1,1)$ entries of both $M$ and $N$ must be positive. In \eqref{eq:rs-sum-B}, this is required for $M_1$ and $M_2$, and then it holds for their sum. These assumptions make the quadratic forms well defined. For the first inequality, subtracting the quadratic form of $M-N$ decreases the expression inside the supremum in \eqref{eq:rs-quadratic-form}. For the second, the expression for $M_1+M_2$ is a sum of two expressions with the same scalar variable $q$; its supremum is at most the sum of their suprema.

For $x>0$ and a direction $(\eta,\xi)$, put
    \begin{equation*}
        \sigma= \left.\frac{\dd}{\dd t}P(y+t\eta,x+t\xi)\right|_{t=0} =P_y(y,x)\eta-f(y,x)\xi.
    \end{equation*}
Thus $\sigma$ is the derivative of a tail probability, rather than a new spatial coordinate. A direct substitution gives
    \begin{equation}
        \label{eq:rs-directional-identity}
        D^2P(y,x)[(\eta,\xi),(\eta,\xi)] =\mathcal B_{\mathcal M(y,x)}(\eta,\sigma).
    \end{equation}
The identity follows by expanding the right-hand side and using $P_x=-f$, $P_{yx}=-f_y$, and $P_{xx}=-f_x$. For $x<0$, we define $\mathcal B_x$ by the algebraic expression $(\sigma-k\eta)^2/J-\ell\eta^2$, with the entries evaluated at $(y,x)$. We do not use the supremum representation in this case, since $J<0$. The relation $P(y,-x)=1-P(y,x)$ gives
    \begin{equation}
        \label{eq:rs-negative-B}
        \mathcal B_{-x}(\eta,\sigma) =-\mathcal B_x(\eta,-\sigma),\qquad x>0,
    \end{equation}
where $\mathcal B_x=\mathcal B_{\mathcal M(y,x)}$ for $x>0$. Indeed, $P(y,-x)=1-P(y,x)$, and all three entries in \eqref{eq:rs-tail-matrix} change sign under $x\mapsto-x$.

Set $u=-a$, $v=b$, and $c=g(y,a,b)$. This is only a linear change of the endpoint variables: it lets us write the defining relation for $g$ as a sum of two tails. The condition $a<b$ becomes $u+v>0$, and we have
    \begin{equation}
        \label{eq:rs-tail-addition}
        P(y,c)=P(y,u)+P(y,v).
    \end{equation}
To prove joint concavity, it is enough to check the second derivative along every straight line in the domain. Fix a point $(y,u,v)$ and constants $\eta,\xi_u,\xi_v$, and set
    \begin{equation*}
        y(t)=y+t\eta,\qquad u(t)=u+t\xi_u,\qquad v(t)=v+t\xi_v, \qquad c(t)=g(y(t),-u(t),v(t)).
    \end{equation*}
For small $|t|$, we still have $u(t)+v(t)>0$. The first three functions are affine, so their second derivatives vanish. All derivatives below are evaluated at $t=0$. Write
    \begin{align*}
        \sigma_u&=\left.\frac{\dd}{\dd t}P(y(t),u(t))\right|_{t=0} = P_y(y,u)\eta-f(y,u)\xi_u,\\
        \sigma_v&=\left.\frac{\dd}{\dd t}P(y(t),v(t))\right|_{t=0} = P_y(y,v)\eta-f(y,v)\xi_v,\\
        \sigma_c&=P_y(y,c)\eta-f(y,c)\dot c.
    \end{align*}
Differentiating \eqref{eq:rs-tail-addition} once gives $\sigma_c=\sigma_u+\sigma_v$. Differentiating twice and using \eqref{eq:rs-directional-identity} gives
    \begin{equation}
        \label{eq:rs-second-directional}
        f(y,c)\ddot c =\mathcal B_c(\eta,\sigma_u+\sigma_v) -\mathcal B_u(\eta,\sigma_u) -\mathcal B_v(\eta,\sigma_v).
    \end{equation}
We prove that this expression is nonpositive. Initially suppose that $u,v,c$ are all nonzero.

\smallskip
\noindent\textit{Case 1: $u,v,c>0$.}
Here the tail probability at $c$ is the sum of the two other positive tail probabilities and is less than $1/2$. Equations \eqref{eq:rs-superadditivity}, \eqref{eq:rs-order-B}, and \eqref{eq:rs-sum-B} show that the right-hand side of \eqref{eq:rs-second-directional} is nonpositive.

\smallskip

\noindent\textit{Case 2: $u,v>0$ and $c<0$.} Let $x_1=u$, $x_2=v$, $x_3=-c$, and set $\sigma_1=\sigma_u$, $\sigma_2=\sigma_v$, $\sigma_3=-\sigma_u-\sigma_v$. By \eqref{eq:rs-negative-B}, it suffices to prove
    \begin{equation*}
        \sum_{i=1}^3\mathcal B_{x_i}(\eta,\sigma_i)\ge0.
    \end{equation*}
Let $J_i,k_i,\ell_i$ be the entries of $\mathcal M(y,x_i)$, and set $J=\sum_{i=1}^3J_i>0$, $k=\sum_{i=1}^3k_i$, and $\ell=\sum_{i=1}^3\ell_i$. The weights in the Cauchy--Schwarz inequality are $J_i>0$:
    \begin{equation*}
        \left(\sum_{i=1}^3(\sigma_i-k_i\eta)\right)^2 = \left(\sum_{i=1}^3 \sqrt{J_i}\,\frac{\sigma_i-k_i\eta}{\sqrt{J_i}}\right)^2 \le \left(\sum_{i=1}^3J_i\right) \left(\sum_{i=1}^3\frac{(\sigma_i-k_i\eta)^2}{J_i}\right).
    \end{equation*}
Since $\sum_{i=1}^3\sigma_i=0$, it follows that
    \begin{align*}
        \sum_{i=1}^3\mathcal B_{x_i}(\eta,\sigma_i) &=\sum_{i=1}^3\frac{(\sigma_i-k_i\eta)^2}{J_i} -\ell\eta^2\\
        &\ge\left(\frac{k^2}{J}-\ell\right)\eta^2.
    \end{align*}
On the other hand, summing \eqref{eq:rs-common-direction} gives
    \begin{equation*}
        J\kappa^2-2k\kappa+\ell\le0.
    \end{equation*}
Consequently
$\ell\le2k\kappa-J\kappa^2\le k^2/J$, proving the claim.

\smallskip

\noindent\textit{Case 3: one of $u,v$ is negative.} By symmetry between $u$ and $v$, assume $u<0<v$. The condition $u+v>0$ implies $v>-u$, and \eqref{eq:rs-tail-addition} gives $c<0$ and
    \begin{equation*}
        P(y,-u)=P(y,v)+P(y,-c)<1/2.
    \end{equation*}
Apply Case 1 to the three positive points $v,-c,-u$. The corresponding derivatives of their tail probabilities are $\sigma_v$, $-\sigma_u-\sigma_v$, and $-\sigma_u$, respectively. Now \eqref{eq:rs-negative-B} gives the required sign in \eqref{eq:rs-second-directional}.

These exhaust the possibilities when $u+v>0$. Points with a zero endpoint or $c=0$ follow by continuity of the Hessian of the implicit threshold, since its defining derivative $f(y,c)$ is strictly positive. Every second directional derivative is nonpositive, proving joint concavity on the convex domain $a<b$.
\end{proof}

\subsection{The profile conditions imply joint concavity}
We now apply Lemma~\ref{lem:rs-matrix-tail} to the normalized sections of the radial density. The assumptions on $w$ will give the required matrix monotonicity and hence the joint concavity in \eqref{eq:JC}.
\begin{lemma}[Joint concavity of the section map]
    \label{lem:profile-JC}
    Under the hypotheses of Theorem~\ref{thm:main}, the function $g$ defined in \eqref{eq:threshold} is jointly concave in the position $y$ and the endpoints $a,b$, as stated in \eqref{eq:JC}.
\end{lemma}

\begin{proof}[Proof of Lemma~\ref{lem:profile-JC}]
    We first take a scalar parameter $y\in\R$ and use the original profile $\varphi$. Set
        \begin{equation*}
            H(y,x)=\varphi(y^2+x^2),\qquad T(y,x)=\int_x^\infty H(y,s)\,ds,\qquad m(y)=\int_{\mathbb R}H(y,s)\,ds.
        \end{equation*}
    The function $x\mapsto H(y,x)$ is not normalized. The functions used in Lemma~\ref{lem:rs-matrix-tail} are
        \begin{equation*}
            f(y,x)=\frac{H(y,x)}{m(y)},\qquad P(y,x)=\frac{T(y,x)}{m(y)}.
        \end{equation*}
    Thus $T(y,0)=m(y)/2$, while $P(y,0)=1/2$. The equation defining $x_y(p)$ can be written as
        \begin{equation*}
            T(y,x_y(p))=p\,m(y),\qquad 0<p<1/2.
        \end{equation*}
    The dependence of $m(y)$ on $y$ must be included when taking derivatives. For example,
        \begin{equation*}
            P_y(y,x)=\frac{T_y(y,x)}{m(y)} -\frac{m'(y)}{m(y)}P(y,x).
        \end{equation*}
    The full matrix calculation is given in \eqref{eq:rs-normalization} below. Both the factor $Z_\varphi^{-1}$ and division by $m(y)$ cancel from the equation defining $g$.

    We first justify the differentiations and limits used below. The assumptions on $w$ give $w\ge w(0)>0$. Also $w(s)/s$ is nonincreasing for $s>0$, so $w(s)=O(1+s)$ and $w'$ is bounded. The identities
        \begin{equation*}
            \varphi'=-\frac{\varphi}{2w},\qquad \varphi''=\frac{1+2w'}{4w^2}\,\varphi
        \end{equation*}
    therefore imply $|\varphi'|+|\varphi''|\le C\varphi$. This justifies differentiating $T$ and $m$ twice in $y$, locally uniformly in $y$, by domination with a constant times $\varphi(s^2)$. The latter is integrable on the line by the radial integrability assumption. Since $\varphi$ is decreasing, the same assumption yields
        \begin{equation*}
            r^n\varphi(r^2)\longrightarrow0\qquad(r\to\infty).
        \end{equation*}
    For example, bound $r^n\varphi(r^2)$ by a constant times $\int_{r/2}^r t^{n-1}\varphi(t^2)\dd t$. These observations justify all tail limits below.

    For $x>0$, define the matrix formed from $H$ and $T$:
        \begin{equation*}
            M(y,x)= \begin{pmatrix} -H^2/H_x&T_y+HH_y/H_x\\ T_y+HH_y/H_x&-T_{yy}-H_y^2/H_x \end{pmatrix}.
        \end{equation*}
    Put $z=y^2+x^2$ and $v=x^2$. Since
        \begin{equation*}
            H_x=-\frac{xH}{w(z)},\qquad H_y=-\frac{yH}{w(z)},
        \end{equation*}
    the entries of $M$ are
        \begin{equation*}
            J=\frac{w(z)H}{x},\qquad k=T_y+\frac{yH}{x},\qquad \ell=-T_{yy}+\frac{y^2H}{w(z)x}.
        \end{equation*}
    They tend to zero as $x\to\infty$. In particular, $J=O(x\varphi(y^2+x^2))\to0$ because $n\ge2$. Direct differentiation gives
        \begin{equation}
            \label{eq:rs-kernel-derivative}
            \partial_x M(y,x)=-H(y,x)\,R(y,x), \quad R(y,x)= \begin{pmatrix} 1-2w'(z)+w(z)/v&y/v\\ y/v&z/(w(z)v) \end{pmatrix}.
        \end{equation}
    For example,
        \begin{equation*}
            J_x=-H\left(1-2w'+\frac w{x^2}\right),\quad k_x=-\frac{yH}{x^2},\quad \ell_x=-\frac{zH}{wx^2}.
        \end{equation*}

    We now show that $R(y,x)$ is nonincreasing in $x$ in the matrix order defined above. For fixed $y$, differentiation with respect to $v=x^2$ gives
        \begin{equation*}
            -v^2\partial_v R = N:={\begin{pmatrix} w-vw'+2v^2w''&y\\ y&\displaystyle\frac{y^2}{w}+\frac{zvw'}{w^2} \end{pmatrix}},
        \end{equation*}
    where $w,w',w''$ are evaluated at $z=y^2+v$. Its first diagonal entry is
        \begin{equation*}
            N_{11}=(w-zw')+y^2w'+2v^2w''\ge0,
        \end{equation*}
    and its second diagonal entry is nonnegative. Its determinant is exactly
        \begin{equation*}
            \det N = \frac{v^2}{w^2} \left\{w'(w-zw') +2w''\bigl(y^2w+zvw'\bigr)\right\}\ge0.
        \end{equation*}
    Thus $N\succeq0$, and $x\mapsto R(y,x)$ is nonincreasing for $x>0$. Integrating \eqref{eq:rs-kernel-derivative} from $x$ to infinity gives
        \begin{equation*}
            M(y,x)=\int_x^\infty H(y,s)R(y,s)\,ds \preceq T(y,x)R(y,x).
        \end{equation*}
    Consequently
        \begin{equation}
            \label{eq:rs-M-over-T}
            \partial_x\left(\frac{M}{T}\right) =\frac{H}{T^2}\bigl(M-TR\bigr)\preceq0.
        \end{equation}
    The integrals are legitimate: for fixed $y$, the entries of $R(y,s)$ are bounded for $s\ge1$, and $H(y,\cdot)$ is integrable.

    We now return to the normalized functions $f=H/m$ and $P=T/m$ defined at the start of this proof. Let $\mathcal M$ be the matrix \eqref{eq:rs-tail-matrix} formed from these functions. Writing $d=m'/m$, a direct differentiation gives
        \begin{equation}
            \label{eq:rs-normalization}
            \frac{\mathcal M}{P} = A^T\frac MT A+C, \quad A=\begin{pmatrix}1&d\\0&1\end{pmatrix}, \quad C=\begin{pmatrix}0&-d\\-d&m''/m-2d^2\end{pmatrix}.
        \end{equation}
    The matrices $A,C$ depend only on $y$, not on $x$. Equation \eqref{eq:rs-M-over-T} implies that $x\mapsto\mathcal M(y,x)/P(y,x)$ is nonincreasing. Since $x\mapsto P(y,x)$ strictly decreases from $1/2$ to $0$ on $(0,\infty)$, this is precisely \eqref{eq:rs-matrix-monotonicity}. The normalized sections meet the regularity assumptions of Lemma~\ref{lem:rs-matrix-tail}; in particular, $f_{xx}(y,0)=2\varphi'(y^2)/m(y)<0$. The lemma proves joint concavity of the threshold in the three scalar variables $(y,a,b)$.

    It remains to pass from a scalar $y$ to a vector $Y\in\R^{n-1}$. We use the same composition argument as at the end of the proof of \cite[Theorem~3.4]{BZ}. Denote the scalar function just obtained by $G(y,a,b)$. For each fixed $a<b$, the preceding argument proves that $y\mapsto G(y,a,b)$ is concave on the whole real line. In particular, this concavity applies to intervals that contain zero. Since $H(-y,x)=H(y,x)$, we also have $G(-y,a,b)=G(y,a,b)$. For fixed $a,b$ and $0\le r\le s$ with $s>0$, write
        \begin{equation*}
            r=\theta s+(1-\theta)(-s), \qquad \theta=\frac{1+r/s}{2}.
        \end{equation*}
    Concavity in $y$ gives
        \begin{equation*}
            G(r,a,b)\ge \theta G(s,a,b)+(1-\theta)G(-s,a,b)=G(s,a,b).
        \end{equation*}
    Thus $G$ is nonincreasing in its first variable on $[0,\infty)$. The actual section at $Y$ is $\varphi(|Y|^2+x^2)$, so its endpoint is $G(|Y|,a,b)$. Put $a_\lambda=(1-\lambda)a_0+\lambda a_1$ and $b_\lambda=(1-\lambda)b_0+\lambda b_1$. The triangle inequality, the monotonicity just proved, and joint concavity of $G$ give
        \begin{align*}
            &G(|(1-\lambda)Y_0+\lambda Y_1|,a_\lambda,b_\lambda)\\
            &\quad\ge G((1-\lambda)|Y_0|+\lambda|Y_1|,a_\lambda,b_\lambda)\\
            &\quad\ge (1-\lambda)G(|Y_0|,a_0,b_0) +\lambda G(|Y_1|,a_1,b_1).
        \end{align*}
    This proves \eqref{eq:JC} in dimension $n$.
\end{proof}

\section{Symmetrization and proof of the main theorem}
\label{sec:geometry}
We now use concavity of the section map to prove \eqref{eq:E} for convex sets. The main step is to find a common sequence of symmetrizations under which the sets converge to parallel half-spaces while keeping their measures.

\subsection{Regularity of sections}

Before defining the symmetrization, we establish continuity of the line integrals and upper semicontinuity of the functions $y\mapsto m_K(y)$ defined below. These properties will ensure that the symmetrized sets are closed.

\begin{lemma}
    \label{lem:regularity}
    Assume the hypotheses of Theorem~\ref{thm:main}. Integrals of $\rho$ over lines and hyperplanes are finite and depend continuously on their position. The functions $(y,s)\mapsto H_y(s)$ and $(y,a,b)\mapsto g(y,a,b)$ are continuous. The distribution function $F$ has a continuous, strictly positive, finite density. For every closed set $K$ and every fixed unit vector $v$, the function
        \begin{equation*}
            m_K(y)=\int_{\R}\ind_K(y+sv)\rho(y+sv)\dd s, \qquad y\in v^\perp,
        \end{equation*}
    is upper semicontinuous.
\end{lemma}

\begin{proof}
    For this proof, write $\vartheta(r)=Z_\varphi^{-1}\varphi(r^2)$, so that $\rho(x)=\vartheta(|x|)$. Normalization gives
        \begin{equation}
            \label{eq:radial-integral}
            \int_0^\infty r^{n-1}\vartheta(r)\dd r<\infty.
        \end{equation}
    Let $1\le k\le n-1$, $|y|\le R$, and $T>\max(2R,1)$. Polar coordinates in $\R^k$, followed by $u=\sqrt{|y|^2+r^2}$, give
        \begin{align}
            &\int_{|z|>T}\vartheta\!\left(\sqrt{|y|^2+|z|^2}\right)\dd z\notag\\
            &\quad=|\Sphere^{k-1}|\int_{\sqrt{|y|^2+T^2}}^\infty u(u^2-|y|^2)^{(k-2)/2}\vartheta(u)\dd u\notag\\
            &\quad\le C_k\int_T^\infty u^{k-1}\vartheta(u)\dd u \longrightarrow0\qquad(T\to\infty), \label{eq:uniform-section-tails}
        \end{align}
    uniformly for $|y|\le R$. For $k=1$, we use $u\ge T>2R\ge2|y|$, which gives
        \begin{equation*}
            \frac{u}{\sqrt{u^2-|y|^2}} = \frac1{\sqrt{1-|y|^2/u^2}} \le\frac1{\sqrt{1-1/4}} = \frac2{\sqrt3}.
        \end{equation*}
    Since $|\Sphere^0|=2$, we can take $C_1=4/\sqrt3$. For $k\ge2$, the bound $(u^2-|y|^2)^{(k-2)/2}\le u^{k-2}$ gives $C_k=|\Sphere^{k-1}|$. Since $u\ge1$ and $k\le n-1$, the right-hand side of \eqref{eq:uniform-section-tails} is finite by \eqref{eq:radial-integral}. On bounded sets, continuity of $\vartheta$ gives a uniform bound. Splitting the integrals into a bounded part and a tail proves continuity of the line integrals and of the hyperplane integrals defining the marginal density.

    For finite $a<b$, the right-hand side of \eqref{eq:threshold} lies strictly between $0$ and $M(y)$. The function $H_y$ is continuous and strictly increasing from $0$ to $M(y)$. Existence and uniqueness of $g$ follow. Joint continuity follows by comparison with $H_y(g-\varepsilon)$ and $H_y(g+\varepsilon)$.

    If $y_j\to y$ and $K$ is closed, then
        \begin{equation*}
            \limsup_j\ind_K(y_j+sv)\le\ind_K(y+sv) \qquad(s\in\R).
        \end{equation*}
    The reverse Fatou inequality on $|s|\le T$, followed by the uniform tail estimate \eqref{eq:uniform-section-tails}, yields $\limsup_jm_K(y_j)\le m_K(y)$.
\end{proof}

\subsection{Geometry, measure preservation and contraction}

Fix a unit vector $v$. For a closed set $K$, its section and the corresponding integral from Lemma~\ref{lem:regularity} are given by
    \begin{equation*}
        K_y=\{s:y+sv\in K\},\qquad m_K(y)=\int_{K_y}\rho(y+sv)\dd s.
    \end{equation*}
The conditional symmetrization of $K$ by left half-lines is
    \begin{equation}
        \label{eq:symmetrization}
        S_vK=\{y+sv:H_y(s)\le m_K(y)\}.
    \end{equation}
If $m_K(y)=0$, the symmetrized section is empty. If $m_K(y)=M(y)$, it is the whole line. If $0<m_K(y)<M(y)$, it is a closed left half-line.

We now use \eqref{eq:JC} to show that symmetrization sends convex sets to convex sets. The same argument gives \eqref{eq:triple-inclusion}, which will allow us to carry a Minkowski inclusion through the whole sequence of symmetrizations.

\begin{lemma}
    \label{lem:basic-symmetrization}
    Assume the hypotheses of Theorem~\ref{thm:main}. For a closed set $K$, the set $S_vK$ is closed, $\mu(S_vK)=\mu(K)$, and $S_v$ is monotone under inclusion. For closed sets $K,L$,
        \begin{equation}
            \label{eq:L1-contraction}
            \mu(S_vK\mathbin{\triangle}S_vL) \le\mu(K\mathbin{\triangle}L).
        \end{equation}
    Since \eqref{eq:JC} holds, $S_vK$ is convex for every closed convex set $K$. Moreover, for closed convex sets,
        \begin{equation}
            \label{eq:triple-inclusion}
            (1-\lambda)K_0+\lambda K_1\subseteq K_2 \quad\Longrightarrow\quad (1-\lambda)S_vK_0+\lambda S_vK_1\subseteq S_vK_2.
        \end{equation}
\end{lemma}

\begin{proof}
    The function $(y,s)\mapsto H_y(s)-m_K(y)$ is lower semicontinuous, so its nonpositive sublevel set \eqref{eq:symmetrization} is closed. By construction, the integral of $\rho$ over each section is unchanged. Tonelli's theorem gives $\mu(S_vK)=\mu(K)$. If $K\subseteq L$, then $m_K(y)\le m_L(y)$ for every $y$, and hence $S_vK\subseteq S_vL$.

    On each line, the two symmetrized sections are nested. The integral of the density over their symmetric difference is $|m_K(y)-m_L(y)|$. Consequently,
        \begin{align*}
            \mu(S_vK\mathbin{\triangle}S_vL) &=\int_{v^\perp}|m_K(y)-m_L(y)|\dd y\\
            &\le\int_{v^\perp}\int_\R |\ind_K(y+sv)-\ind_L(y+sv)|\rho(y+sv)\dd s\dd y,
        \end{align*}
    which proves \eqref{eq:L1-contraction}. 
    For \eqref{eq:triple-inclusion}, take $y_i+s_iv\in S_vK_i$, $i=0,1$, and $\varepsilon>0$. Recall that
        \begin{equation*}
            (K_i)_{y_i} = \{t\in\R:y_i+tv\in K_i\}, \qquad i=0,1.
        \end{equation*}
    These sections are closed intervals, possibly unbounded. Since $H_{y_i}(s_i-\varepsilon)<m_{K_i}(y_i)$, we can choose finite nondegenerate intervals $[a_i,b_i]\subseteq(K_i)_{y_i}$ 
    % \tz{<--Czy było już wcześniej takie oznaczenie?}  
    such that
        \begin{equation*}
            g(y_i,a_i,b_i)>s_i-\varepsilon.
        \end{equation*}
    This also applies to unbounded and full sections by truncation. Use a subscript $\lambda$ for the corresponding affine combinations. The assumed inclusion gives $[a_\lambda,b_\lambda]\subseteq(K_2)_{y_\lambda}$. Joint concavity implies
        \begin{equation*}
            g(y_\lambda,a_\lambda,b_\lambda) \ge(1-\lambda)g(y_0,a_0,b_0)+\lambda g(y_1,a_1,b_1) >s_\lambda-\varepsilon.
        \end{equation*}
    Thus $y_\lambda+(s_\lambda-\varepsilon)v\in S_vK_2$. Letting $\varepsilon\downarrow0$ and using closedness proves \eqref{eq:triple-inclusion}. Taking $K_0=K_1=K_2$ proves convexity.
\end{proof}

\begin{remark}[Meaning of the analytic hypothesis]
    \label{rem:equivalence}
    Under the assumptions of Theorem~\ref{thm:main}, condition \eqref{eq:JC} is also necessary for $S_vK$ to be convex whenever $K$ is a closed convex set. Indeed, for $y_0\ne y_1$, take the convex hull of $\{y_0+sv:a_0\le s\le b_0\}$ and $\{y_1+sv:a_1\le s\le b_1\}$. Its section at $y_\lambda$ is exactly $[a_\lambda,b_\lambda]$. Convexity of the symmetrized set gives \eqref{eq:JC-expanded}. The case $y_0=y_1$ follows by perturbation and continuity. When $n>2$, this argument uses lower-dimensional convex sets. These are included in the assertion that $S_vK$ is convex for every closed convex set $K$. This equivalence concerns symmetrization, not necessity of \eqref{eq:JC} for \eqref{eq:E}.
\end{remark}

\subsection{Compactness at fixed positive measure}

To study repeated symmetrizations, we need convergent subsequences of closed convex sets with a fixed measure in $(0,1)$. The next lemma also shows that Minkowski inclusions pass to simultaneous limits of the sets.

\begin{lemma}
    \label{lem:compactness}
    For $0<p<1$, let $\mathcal C_p$ be the family of closed convex sets $K$ with $\mu(K)=p$. Every sequence in $\mathcal C_p$ has a subsequence, still denoted by $(K_j)$, and a closed convex set $K$ such that
        \begin{equation*}
            \dist(\cdot,K_j)\longrightarrow\dist(\cdot,K) \qquad\text{locally uniformly on }\R^n.
        \end{equation*}
    For this limit set $K$, we also have
        \begin{equation}
            \label{eq:measure-convergence}
            \mu(K_j\mathbin{\triangle}K)\longrightarrow0, \qquad \mu(K)=p.
        \end{equation}
    Moreover, the relation $(1-\lambda)A+\lambda B\subseteq C$ is closed under simultaneous convergence of the three distance functions.
\end{lemma}

\begin{proof}
    Write $B_R=RB^n=\{x\in\R^n:|x|\le R\}$ and choose $R$ such that $\mu(\R^n\setminus B_R)<p$. Every $K_j\in\mathcal C_p$ intersects $B_R$. Hence $d_j=\dist(\cdot,K_j)$ are convex, nonnegative, $1$-Lipschitz functions with $d_j(0)\le R$. The Arzel\`a-Ascoli theorem and a diagonal argument yield a locally uniform limit $d$. The set $K=\{d=0\}$ is closed, convex, and nonempty: choose $z_j\in K_j\cap B_R$ and pass to a convergent subsequence. The Lipschitz property gives $d(x)\le\dist(x,K)$. Conversely, for fixed $x$, nearest points $z_j\in K_j$ remain bounded because $|x-z_j|=d_j(x)$. Any limit point $z$ belongs to $K$ and satisfies $|x-z|=d(x)$. Thus $d=\dist(\cdot,K)$. Every $z\in K$ can consequently be approximated by $z_j\in K_j$.

    If $x\notin K$, then $x\notin K_j$ eventually. If $K$ had empty interior, it would lie in a proper affine subspace, and dominated convergence would give $\mu(K_j)\to0$, a contradiction. Thus $K$ has nonempty interior. For $x\in\Int K$, choose the vertices of a small full-dimensional simplex in $K$ surrounding $x$, and approximate them by points of $K_j$. For all sufficiently large $j$, $x$ lies in the simplex of approximating vertices, and hence in $K_j$. Therefore $\ind_{K_j}\to\ind_K$ outside $\partial K$. The boundary of a convex set with nonempty interior has Lebesgue measure zero. Since $\mu$ has a density, it also has $\mu$-measure zero. Dominated convergence proves \eqref{eq:measure-convergence}.

    Finally, if $(1-\lambda)A_j+\lambda B_j\subseteq C_j$, approximate $a\in A$ and $b\in B$ by $a_j\in A_j$ and $b_j\in B_j$. Then $(1-\lambda)a_j+\lambda b_j\in C_j$. Taking the limit in the distance to $C_j$ proves $(1-\lambda)a+\lambda b\in C$.
\end{proof}

\subsection{Strict improvement of a bounded functional}

Fix $u\in\Sphere^{n-1}$ and define
    \begin{equation*}
        J_u(K)=-\int_K\arctan(u\cdot x)\dd\mu(x).
    \end{equation*}
Then $|J_u(K)|\le\pi/2$ and
    \begin{equation}
        \label{eq:J-Lipschitz}
        |J_u(K)-J_u(L)|\le\frac\pi2\mu(K\mathbin{\triangle}L).
    \end{equation}

    The function $\arctan$ is bounded, so $J_u$ is finite without any moment assumption on $\mu$. For $u\cdot v>0$, the next lemma shows that $J_u(S_vK)>J_u(K)$ whenever $S_vK\ne K$.
\begin{lemma}
    \label{lem:strict-improvement}
    Assume the hypotheses of Theorem~\ref{thm:main}. For every closed convex set $K$ with $0<\mu(K)<1$ and every pair of unit vectors $u,v$ satisfying $u\cdot v>0$, we have
        \begin{equation}
            \label{eq:strict-improvement}
            J_u(S_vK)\ge J_u(K),
        \end{equation}
    with equality if and only if $S_vK=K$.
\end{lemma}

\begin{proof}
    On the line above $y$, the function $\theta_y(s)=-\arctan(u\cdot y+s\,u\cdot v)$ is strictly decreasing. Suppose first that $0<m_K(y)<M(y)$. Put $I=K_y$ and $L=(-\infty,h]$, where $H_y(h)=m_K(y)$. Then
        \begin{align*}
            &\int_L\theta_y(s)\rho(y+sv)\dd s -\int_I\theta_y(s)\rho(y+sv)\dd s\\
            &\quad=\int_\R(\ind_L-\ind_I) (\theta_y(s)-\theta_y(h))\rho(y+sv)\dd s\ge0.
        \end{align*}
    Both factors have the same sign on each side of $h$. Since $\rho>0$ and $\theta_y$ is strictly decreasing, equality holds exactly when the symmetric difference of $I$ and $L$ has measure zero for the density $s\mapsto\rho(y+sv)$. The conclusion also holds when $m_K(y)=0$ or $m_K(y)=M(y)$. Integrating over $y$ gives \eqref{eq:strict-improvement}, with equality precisely when $\mu(S_vK\mathbin{\triangle}K)=0$.

    By Lemma~\ref{lem:basic-symmetrization}, both $K$ and $S_vK$ are closed and convex, and their measures are positive. Each of these sets is therefore the closure of its interior. If they were distinct, their symmetric difference would contain a nonempty open set. Since $\rho>0$, this open set would have positive $\mu$-measure. Thus $\mu(S_vK\mathbin{\triangle}K)=0$ implies $S_vK=K$.
\end{proof}

\subsection{Choosing a common sequence of symmetrizations}

To preserve a Minkowski inclusion at every step, we must apply the same symmetrization to all the sets involved. The following construction uses the sum of their functionals $J_u$ to choose a common sequence of directions and obtain parallel half-spaces as limits.

\begin{proposition}[Simultaneous convergence of conditional symmetrizations]
    \label{prop:simultaneous-convergence}
    Assume the hypotheses of Theorem~\ref{thm:main}. For every unit vector $u$ and every finite family $K_1,\ldots,K_m$ of closed convex sets with $0<\mu(K_i)<1$, there exists a common sequence of unit directions $w_1,w_2,\ldots$, with $u\cdot w_k>0$, such that
        \begin{equation*}
            K_i^{(0)}=K_i,\qquad K_i^{(k)}=S_{w_k}K_i^{(k-1)}
        \end{equation*}
    satisfy
        \begin{align*}
                \mu(K_i^{(k)}\mathbin{\triangle}H_i)&\longrightarrow0,\\
                \dist(\cdot,K_i^{(k)})&\longrightarrow\dist(\cdot,H_i) \quad\text{locally uniformly},
        \end{align*}
    where
        \begin{equation}
            \label{eq:limit-halfspaces}
            H_i=\{x:u\cdot x\le F^{-1}(\mu(K_i))\}.
        \end{equation}
    The sequence may depend on the initial family and on $u$.
\end{proposition}

\begin{proof}
    Choose a countable dense family $(v_j)_{j\ge1}$ in
        \begin{equation*}
            \mathcal H_u=\{v\in\Sphere^{n-1}:u\cdot v>0\}.
        \end{equation*}
    For a tuple $\mathcal K=(K_1,\ldots,K_m)$, set
        \begin{equation}
            \label{eq:gain}
            \Delta_j(\mathcal K) =\sum_{i=1}^m\bigl(J_u(S_{v_j}K_i)-J_u(K_i)\bigr)\ge0.
        \end{equation}
    Beginning with $\mathcal K^{(0)}$, at step $k\ge1$ choose
        \begin{equation*}
            j_k\in\argmax_{1\le j\le k}\Delta_j(\mathcal K^{(k-1)}), \qquad w_k=v_{j_k},\qquad K_i^{(k)}=S_{w_k}K_i^{(k-1)}.
        \end{equation*}
    The maximum is over a finite set; no continuity in the direction parameter is required. The same direction is applied to all $m$ sets. By Lemma~\ref{lem:basic-symmetrization}, all sets remain closed and convex, and their measures remain unchanged. The sequence $\mathcal J_k=\sum_iJ_u(K_i^{(k)})$ is nondecreasing and bounded. Its increments
        \begin{equation*}
            \delta_k=\mathcal J_k-\mathcal J_{k-1} =\Delta_{j_k}(\mathcal K^{(k-1)})
        \end{equation*}
    tend to zero. For every fixed $j$, and every $k\ge j$,
        \begin{equation}
            \label{eq:fixed-gain-zero}
            0\le\Delta_j(\mathcal K^{(k-1)})\le\delta_k\longrightarrow0.
        \end{equation}

    Lemma~\ref{lem:compactness}, applied to the finite tuple, gives convergent subsequences in the distance-function topology and in $L^1(\mu)$. Let $\mathcal K^*=(K_1^*,\ldots,K_m^*)$ be any such limit. By \eqref{eq:L1-contraction} and \eqref{eq:J-Lipschitz}, each fixed $\Delta_j$ is continuous in $L^1(\mu)$; indeed its change is at most $\pi\sum_i\mu(K_i\mathbin{\triangle}L_i)$. By \eqref{eq:fixed-gain-zero}, we have $\Delta_j(\mathcal K^*)=0$. Every summand in \eqref{eq:gain} is nonnegative, so Lemma~\ref{lem:strict-improvement} implies
        \begin{equation}
            \label{eq:fixed-points}
            S_{v_j}K_i^*=K_i^*\qquad\text{for all }i,j.
        \end{equation}

    For each $i$ and $j$, every section of $S_{v_j}K_i^*$ is a left half-line, the whole line, or the empty set. Hence \eqref{eq:fixed-points} gives
        \begin{equation*}
            K_i^*-t v_j\subseteq K_i^*, \qquad t\ge0.
        \end{equation*}
    Since the directions $v_j$ are dense in $\mathcal H_u$ and $K_i^*$ is closed, this inclusion holds for every unit vector $v$ with $u\cdot v\ge0$. In particular, it holds for both signs of every direction in $u^\perp$. Thus $K_i^*$ is invariant under translations in $u^\perp$, and $K_i^*-tu\subseteq K_i^*$ for every $t\ge0$. Since $0<\mu(K_i^*)<1$, these properties imply
        \begin{equation*}
            K_i^*=\{x:u\cdot x\le c_i\}
        \end{equation*} 
    for some finite $c_i$. Finally, $\mu(K_i^*)=\mu(K_i)$ determines $c_i=F^{-1}(\mu(K_i))$ uniquely.

    Thus every subsequential limit is the tuple of half-spaces in \eqref{eq:limit-halfspaces}. By Lemma~\ref{lem:compactness}, every subsequence has a convergent further subsequence. Since the limit is unique, the whole sequence of distance functions converges to the distance functions of these half-spaces. As in the proof of Lemma~\ref{lem:compactness}, this convergence and convexity imply
        \begin{equation*}
            \ind_{K_i^{(k)}}(x)\longrightarrow\ind_{H_i}(x) \qquad\text{for }x\notin\partial H_i.
        \end{equation*}
    The hyperplane $\partial H_i$ has $\mu$-measure zero. Dominated convergence therefore gives
        \begin{equation*}
            \mu(K_i^{(k)}\mathbin{\triangle}H_i) \longrightarrow0.
        \end{equation*}
    This proves both asserted convergences.
\end{proof}

\subsection{Transfer of Minkowski inclusion to the limit}

We now apply the same symmetrizations to the three sets in a Minkowski inclusion and pass to their limiting parallel half-spaces. The formula for the Minkowski combination of two parallel half-spaces then gives \eqref{eq:E}.

\begin{proof}[Proof of Theorem~\ref{thm:main}]
    Lemma~\ref{lem:profile-JC} proves joint concavity of $g$, that is, \eqref{eq:JC}. Taking closures does not change the measures of the convex sets, so it is enough to consider closed convex $K,L$. Put $C=\overline{(1-\lambda)K+\lambda L}$. If $\mu(C)=1$, there is nothing to prove. Otherwise $\mu(K)$, $\mu(L)$, and $\mu(C)$ belong to $(0,1)$. Apply Proposition~\ref{prop:simultaneous-convergence} to the tuple $(K,L,C)$, using the same symmetrization direction at every step. The initial inclusion $(1-\lambda)K+\lambda L\subseteq C$ remains valid throughout the iteration by \eqref{eq:triple-inclusion}. It passes to the limit by Lemma~\ref{lem:compactness}. The limit sets are the parallel half-spaces $H_{c_0},H_{c_1},H_{c_2}$, where
        \begin{equation*}
            c_0=\Psi(\mu(K)),\qquad
            c_1=\Psi(\mu(L)),\qquad
            c_2=\Psi(\mu(C)).
        \end{equation*}
    For $H_c=\{x:u\cdot x\le c\}$, the limiting
    inclusion reads
        \begin{equation*}
            (1-\lambda)H_{c_0}+\lambda H_{c_1} = H_{(1-\lambda)c_0+\lambda c_1} \subseteq H_{c_2}.
        \end{equation*}
    Therefore
        \begin{equation*}
            (1-\lambda)c_0+\lambda c_1\le c_2,
        \end{equation*} 
    which is \eqref{eq:E}.
    % Since $(1-\lambda)H_{c_0}+\lambda H_{c_1} =H_{(1-\lambda)c_0+\lambda c_1}$, their inclusion gives $(1-\lambda)c_0+\lambda c_1\le c_2$. \tz{Czy tu nie ma po prostu równości??} This is \eqref{eq:E}.
\end{proof}

\section{Examples: Gaussian, Student, and a non-Student family}
\label{sec:examples}

We first check the criterion for Gaussian and Student measures. We then give a family for which $w$ is not affine. In each example, \eqref{eq:E} uses the inverse of the one-dimensional marginal distribution function of that measure.

\subsection{Gaussian measures}

\begin{theorem}
    \label{thm:gaussian}
    For every $\beta>0$, the probability measure
        \begin{equation*}
            \dd\gamma_{n,\beta}(x) =\left(\frac{\beta}{\pi}\right)^{n/2}e^{-\beta|x|^2}\dd x
        \end{equation*}
    satisfies \eqref{eq:E}, with $F(t)=\Phi(\sqrt{2\beta}\,t)$, where $\Phi$ is the standard Gaussian distribution function.
\end{theorem}

\begin{proof}
    Here $\varphi(s)=e^{-\beta s}$ and $w(s)\equiv(2\beta)^{-1}$. Therefore $w'=w''=0$ and $w-sw'=w>0$. Theorem~\ref{thm:main} applies.
\end{proof}

\subsection{Multivariate Student measures}
    
\begin{theorem}
    \label{thm:student}
    Let $\tau>0$ and $p>n/2$. Then
        \begin{equation*}
            \dd\mu_{n,p,\tau}(x) = \frac{\Gamma(p)\tau^{p-n/2}}{\pi^{n/2}\Gamma(p-n/2)}(\tau+|x|^2)^{-p}\dd x
        \end{equation*}
    satisfies \eqref{eq:E} with the inverse of its one-dimensional marginal distribution function.
\end{theorem}

\begin{proof}
    Take $\varphi(s)=(\tau+s)^{-p}$. Then
        \begin{equation*}
            w(s)=\frac{\tau+s}{2p},\quad w'(s)=\frac1{2p}>0,\quad w''(s)=0,\quad w(s)-sw'(s)=\frac\tau{2p}>0.
        \end{equation*}
    Euler's Beta integral gives
        \begin{equation*}
            \int_0^\infty r^{n-1}(\tau+r^2)^{-p}\dd r = \frac{\tau^{n/2-p}}2 \frac{\Gamma(n/2)\Gamma(p-n/2)}{\Gamma(p)} <\infty \qquad\left(p>\frac n2\right).
        \end{equation*}
    Theorem~\ref{thm:main} gives the assertion.
\end{proof}

The standard isotropic Cauchy law corresponds to $p=(n+1)/2$ and $\tau=1$ in the Student family.

\begin{corollary}[Cauchy measures and boxes]
    \label{cor:cauchy}
    For $p=(n+1)/2$ and $\tau=1$, the density
        \begin{equation*}
            \dd\mu_{n,1}(x) = \frac{\Gamma((n+1)/2)}{\pi^{(n+1)/2}} (1+|x|^2)^{-(n+1)/2}\dd x
        \end{equation*}
    is the isotropic Cauchy law, with marginal $F_1(t)=\frac12+\frac1\pi\arctan t$ and $F_1^{-1}(p)=-\cot(\pi p)$. It satisfies \eqref{eq:E} for all convex sets with measures in $(0,1)$. In particular, for axis-parallel boxes the function
        \begin{equation*}
            (a_1,b_1,\ldots,a_n,b_n)\longmapsto \cot\!\left(\pi\mu_{n,1} \left(\prod_{j=1}^n(a_j,b_j)\right)\right)
        \end{equation*}
    is jointly convex on the domain $a_j<b_j$, $j=1,\ldots,n$. Here $\prod_{j=1}^n(a_j,b_j)=(a_1,b_1)\times\cdots\times(a_n,b_n)$ denotes the Cartesian product of the intervals.
\end{corollary}

\begin{proof}
    The marginal is standard Cauchy. Theorem~\ref{thm:student} gives the inequality. Minkowski interpolation of boxes interpolates their endpoints; multiplication by $-1$ yields the convexity statement.
\end{proof}

\subsection{Adjacent-power mixtures beyond the Student class}

The next family combines two Student profiles and gives a nonaffine function $w$ when $\varepsilon>0$. We verify the sufficient conditions here; Section~\ref{sec:necessary} will show that the resulting parameter range is exact.

\begin{corollary}
    \label{cor:rs-adjacent-mixtures}
    Fix $n\ge2$, $p>n/2$, $\tau>0$, and
        \begin{equation*}
            0\le\varepsilon\le p\tau.
        \end{equation*}
    The probability density
        \begin{equation}
            \label{eq:rs-mixture-density}
            \rho_{n,p,\tau,\varepsilon}(x) =\frac1{Z_{n,p,\tau,\varepsilon}} \left[(\tau+|x|^2)^{-p}+\varepsilon(\tau+|x|^2)^{-p-1}\right],
        \end{equation}
    where
        \begin{equation*}
            Z_{n,p,\tau,\varepsilon} =\frac{\pi^{n/2}\Gamma(p-n/2)}{\Gamma(p)} \tau^{n/2-p} \left(1+\frac{\varepsilon(p-n/2)}{p\tau}\right),
        \end{equation*}
    satisfies \eqref{eq:JC} and hence \eqref{eq:E}. If $\varepsilon>0$, this density is neither Gaussian nor any individual radial Student density, including rescaled ones.
\end{corollary}

\begin{proof}
    The normalizing constant follows from polar coordinates and Euler's Beta integral. The squared-radius profile
        \begin{equation*}
            \varphi(s)=(\tau+s)^{-p}+\varepsilon(\tau+s)^{-p-1}
        \end{equation*}
    is positive, smooth, strictly decreasing, and integrable in dimension $n$ because $p>n/2$. Its reciprocal logarithmic slope is
        \begin{equation*}
            w(s)=\frac{(\tau+s)(\tau+s+\varepsilon)}{2[p(\tau+s)+(p+1)\varepsilon]}.
        \end{equation*}
    Writing $A=\tau+s$ and $D=pA+(p+1)\varepsilon$, differentiation gives
        \begin{align}
            w'(s)&=\frac{pA^2+2(p+1)\varepsilon A+(p+1)\varepsilon^2} {2D^2}>0,\nonumber\\
            w''(s)&=\frac{(p+1)\varepsilon^2}{D^3}\ge0. \label{eq:rs-mixture-wsecond}
        \end{align}
    The useful division formula is
        \begin{equation*}
            w(s)=\frac{s}{2p} +\frac{p\tau-\varepsilon}{2p^2} +\frac{(p+1)\varepsilon^2}{2p^2D}.
        \end{equation*}
    In particular,
        \begin{equation*}
            \lim_{s\to\infty}\bigl(w(s)-s w'(s)\bigr) =\frac{p\tau-\varepsilon}{2p^2}\ge0.
        \end{equation*}
    Since
        \begin{equation*}
            \frac{d}{ds}(w-sw')=-s w''\le0,
        \end{equation*}
    we have $w-sw'\ge0$ for every $s\ge0$. Theorem~\ref{thm:main} applies.

    For $\varepsilon>0$, equation \eqref{eq:rs-mixture-wsecond} gives $w''(s)>0$ everywhere. For every Gaussian profile $w$ is constant, and for every radial Student profile it is affine. Multiplication of the profile by a normalizing constant leaves $w$ unchanged. Thus \eqref{eq:rs-mixture-density} cannot be a member of either individual family.
\end{proof}

\begin{remark}
    The proof verifies \eqref{eq:JC} directly for \eqref{eq:rs-mixture-density} when $0\le\varepsilon\le p\tau$. It does not use or assert preservation of \eqref{eq:E} under arbitrary mixtures. Failure for $\varepsilon>p\tau$ is proved in Section~\ref{sec:necessary}.
\end{remark}

\subsection{An explicit planar example}

We now take $n=2$, $p=3/2$, and $\tau=\varepsilon=1$ in Corollary~\ref{cor:rs-adjacent-mixtures}. For this choice, both the density and its one-dimensional marginal distribution function have simple formulas.

\begin{proposition}
    \label{prop:rs-explicit-planar}
    The probability density
        \begin{equation*}
            \rho_\dagger(x,y) = \frac{3(2+x^2+y^2)}{8\pi(1+x^2+y^2)^{5/2}}
        \end{equation*}
    is radial, satisfies \eqref{eq:JC}, and is neither Gaussian nor an individual Student density. Its marginal distribution function is
        \begin{equation}
            \label{eq:rs-explicit-marginal}
            F_\dagger(t)=\frac12+\frac1\pi\arctan t + \frac{t}{4\pi(1+t^2)}.
        \end{equation}
    Let $\mu_\dagger$ be the probability measure with density $\rho_\dagger$. Then, for all convex Borel sets $K,L\subseteq\mathbb R^2$ such that $0<\mu_\dagger(K),\mu_\dagger(L)<1$, and every $0<\lambda<1$,
        \begin{equation*}
            F_\dagger^{-1}\!\left(\mu_\dagger((1-\lambda)K+\lambda L)\right) \ge(1-\lambda)F_\dagger^{-1}(\mu_\dagger(K))+\lambda F_\dagger^{-1}(\mu_\dagger(L)).
        \end{equation*}
    If $\mu_\dagger((1-\lambda)K+\lambda L)=1$, the left-hand side is interpreted as $F_\dagger^{-1}(1)=+\infty$.
\end{proposition}

\begin{proof}
    This is Corollary~\ref{cor:rs-adjacent-mixtures} with $n=2$, $p=3/2$, $\tau=1$, and $\varepsilon=1$. For completeness, omitting its irrelevant positive normalizing factor, the squared-radius profile is
        \begin{equation*}
            \varphi(s)=\frac{s+2}{(s+1)^{5/2}}.
        \end{equation*}
    The three conditions can be checked explicitly:
        \begin{align*}
            w(s)&=\frac{(s+1)(s+2)}{3s+8},\\
            w'(s)&=\frac{3s^2+16s+18}{(3s+8)^2}>0,\\
            w''(s)&=\frac{20}{(3s+8)^3}>0,\\
            w(s)-s w'(s)&=\frac{s^2+12s+16}{(3s+8)^2}>0.
        \end{align*}
    Integrating the density in $y$ gives
        \begin{equation*}
            f_\dagger(t)=\frac{3}{4\pi(1+t^2)} + \frac{1}{2\pi(1+t^2)^2},
        \end{equation*}
    whose integral is \eqref{eq:rs-explicit-marginal}. The claimed inequality follows from the main theorem.
\end{proof}

\begin{remark}
    In the notation of the Cauchy--Student mixture
        \begin{equation*}
            \rho_\varepsilon(x,y) = \frac{(1+x^2+y^2)^{-3/2} +\varepsilon(1+x^2+y^2)^{-5/2}} {2\pi(1+\varepsilon/3)},
        \end{equation*}
    Corollary~\ref{cor:rs-adjacent-mixtures} verifies the full hypothesis \eqref{eq:JC} for every $0\le\varepsilon\le3/2$. This is joint concavity in all variables $(y,a,b)$, not merely concavity in $a,b$ for each fixed $y$. The resulting Ehrhard-type inequality uses the inverse of the marginal distribution function of $\rho_\varepsilon$, not the inverse Cauchy distribution function.
\end{remark}

\section{Necessary conditions and sharp obstructions}
\label{sec:necessary}

We now prove the necessary conditions stated in the introduction. We assume neither the conditions on $w$ nor any particular form of the one-dimensional densities. Throughout this part, $\varphi:[0,\infty)\to[0,\infty)$ is continuous, not identically zero, and $0<Z_\varphi<\infty$. Thus the profile may initially have zeros. As in \eqref{eq:R}, $\rho(x)=Z_\varphi^{-1}\varphi(|x|^2)$, and $F$ is the distribution function of the first coordinate. The endpoint values $F^{-1}(0)$ and $F^{-1}(1)$ are the one-sided limits, which may be infinite. We require \eqref{eq:E} when both $K$ and $L$ have measures strictly between $0$ and $1$. In particular, if their Minkowski combination has measure zero, then \eqref{eq:E} fails. We use $B^n=\{x\in\R^n:|x|\le1\}$ for the Euclidean unit ball.

We begin with small balls whose centers lie on the same sphere. Their measures are equal by radial symmetry, so a midpoint comparison in \eqref{eq:E}, followed by shrinking the balls, will force $\varphi$ to be nonincreasing.

\begin{proposition}[Radial monotonicity is necessary]
    \label{prop:necessary-monotonicity}
    If $\mu$ satisfies \eqref{eq:E} for all convex sets $K,L$ with $0<\mu(K),\mu(L)<1$, then \(\varphi\) is nonincreasing on \([0,\infty)\).
\end{proposition}

\begin{proof}
    Fix \(0\le s<t\) with \(\varphi(t)>0\). Since \(n\ge2\), there exist \(x_0,x_1\in\mathbb R^n\) such that
        \begin{equation*}
            |x_0|^2=|x_1|^2=t,\qquad \left|\frac{x_0+x_1}{2}\right|^2=s.
        \end{equation*}
    For \(r>0\) let \(K_i=x_i+rB^n\). By rotational invariance, \(\mu(K_0)=\mu(K_1)=p_r>0\) for all sufficiently small \(r\), while
        \begin{equation*}
            \frac{K_0+K_1}{2}=\frac{x_0+x_1}{2}+rB^n.
        \end{equation*}
    Since $F^{-1}$ is increasing, \eqref{eq:E} gives \(\mu((K_0+K_1)/2)\ge p_r\). Dividing by \(|B^n|r^n\) and letting \(r\downarrow0\) yields \(\varphi(s)\ge\varphi(t)\). If \(\varphi(t)=0\), the same inequality holds trivially. This proves the claim.
\end{proof}

We next perturb a half-space by the graph of the convex function $y\mapsto|y|$, using the monotonicity just proved to justify differentiation under the integral. The convexity property obtained in the next lemma will rule out bounded support in the proof of Theorem~\ref{thm:necessary-full-support}.

\begin{lemma}[Conditional convexity]
    \label{lem:conditional-moment-convexity}
    Assume \eqref{eq:E}, and write a random vector with law \(\mu\) as \(Z=(X,Y)\in\mathbb R\times\mathbb R^{n-1}\). At every \(t\) where the marginal density of \(X\) is positive, set
        \begin{equation*}
            m(t)=\mathbb E\bigl[|Y|\mid X=t\bigr].
        \end{equation*}
    Then $m(-t)=m(t)$, and $m$ is convex on the interior of the support of the marginal distribution. Whenever the denominator is positive, it is given by
        \begin{equation}
            \label{eq:conditional-profile-ratio}
            m(t)= \frac{\displaystyle\int_0^\infty u^{(n-2)/2}\varphi(t^2+u)\,du}{\displaystyle\int_0^\infty u^{(n-3)/2}\varphi(t^2+u)\,du}.
        \end{equation}
\end{lemma}

\begin{proof}
    For \(\varepsilon\ge0\) and \(t\in\mathbb R\), put
        \begin{equation*}
            K_{\varepsilon,t} = \{(x,y)\in\mathbb R\times\mathbb R^{n-1}: x+\varepsilon|y|\le t\}.
        \end{equation*}
    For fixed \(\varepsilon\), these sets are translates of a common convex set. Consequently
        \begin{equation*}
            (1-\lambda)K_{\varepsilon,t_0} +\lambda K_{\varepsilon,t_1} =K_{\varepsilon,(1-\lambda)t_0+\lambda t_1}.
        \end{equation*}
    For each fixed $\varepsilon\ge0$, inequality \eqref{eq:E} implies that $U_\varepsilon(t):=F^{-1}(\mu(K_{\varepsilon,t}))$ is concave on the interval where $0<\mu(K_{\varepsilon,t})<1$. Since \(U_0(t)=t\), differentiation at \(\varepsilon=0+\) gives
        \begin{equation*}
            \left.\partial_\varepsilon U_\varepsilon(t)\right|_{0+} = -\frac{\int_{\mathbb R^{n-1}}|y|\rho(t,y)\,dy}{\int_{\mathbb R^{n-1}}\rho(t,y)\,dy} =-m(t).
        \end{equation*}
    To justify the differentiation, write the difference of the measures as an integral over $t-\varepsilon|y|<x\le t$. Since $\varphi$ is nonincreasing, for every $\varepsilon>0$ its difference quotient is bounded by
        \begin{equation*}
            \frac1\varepsilon \int_{t-\varepsilon|y|}^{t}\rho(x,y)\dd x \le \frac{|y|}{Z_\varphi}\varphi(|y|^2).
        \end{equation*}
    This bound is integrable on $\R^{n-1}$, since
        \begin{equation*}
            \int_{\R^{n-1}}|y|\varphi(|y|^2)\dd y =|\Sphere^{n-2}|\int_0^\infty r^{n-1}\varphi(r^2)\dd r <\infty.
        \end{equation*}
    Dominated convergence therefore gives the displayed derivative. No first-moment assumption on $\mu$ is needed. Taking the first-order limit of the concavity inequality for \(U_\varepsilon\) gives convexity of \(m\). Rotational symmetry gives \(m(t)=m(-t)\). Formula \eqref{eq:conditional-profile-ratio} follows by using polar coordinates in \(Y\) and setting \(u=|y|^2\).
\end{proof}

\begin{proof}[Proof of Theorem~\ref{thm:necessary-full-support}]
    Suppose, to the contrary, that the support of $\mu$ is bounded. Let
        \begin{equation*}
            R=\sup\{|x|:x\in\operatorname{supp}\mu\}<\infty.
        \end{equation*}
    Thus $R$ is the smallest radius of a closed ball centered at the origin that contains the support. By Proposition~\ref{prop:necessary-monotonicity}, $\varphi$ is nonincreasing. The marginal density of $X$ is positive on $(-R,R)$, and Lemma~\ref{lem:conditional-moment-convexity} gives $m(0)>0$. Since $m(-t)=m(t)$ and $m$ is convex,
        \begin{equation*}
            m(0)\le\frac{m(t)+m(-t)}2=m(t).
        \end{equation*}
    Consequently,
        \begin{equation*}
            m(t)\ge m(0)>0\qquad(0\le t<R).
        \end{equation*}
    On the other hand, conditioned on \(X=t\), the support restriction gives \(|Y|\le\sqrt{R^2-t^2}\). Hence \(m(t)\le\sqrt{R^2-t^2}\to0\) as \(t\uparrow R\), a contradiction. Thus the support is not bounded. Proposition~\ref{prop:necessary-monotonicity} shows that \(\varphi\) is nonincreasing. If \(\varphi(s_0)=0\) at any finite \(s_0\), nonnegativity and monotonicity would give \(\varphi(s)=0\) for every \(s\ge s_0\), contradicting the unboundedness of the support. This proves all three assertions.
\end{proof}

We now assume the power-law marginal tail in Theorem~\ref{prop:necessary-small-balls} and apply \eqref{eq:E} to balls whose radii tend to zero. This gives a convexity condition on the density, which will then be used to restrict the second term of the radial tail expansion.

\begin{proof}[Proof of Theorem~\ref{prop:necessary-small-balls}]
    Let $v_n=|B^n|$ be the volume of the Euclidean unit ball. For fixed $x\in\R^n$ and $r>0$, as $\varepsilon\downarrow0$, we have
        \begin{equation*}
            \mu(x+\varepsilon rB^n) = v_n\varepsilon^n r^n\rho(x)(1+o(1)).
        \end{equation*}
    Symmetry and \eqref{eq:necessary-marginal-tail} give \(F^{-1}(p)\sim-(A/p)^{1/\eta}\) as \(p\downarrow0\). Minkowski interpolation of balls interpolates both their centers and radii. Thus \eqref{eq:E} gives concavity of $(x,r)\mapsto F^{-1}(\mu(x+\varepsilon rB^n))$. Multiplication by the negative number $-\varepsilon^{n/\eta}$ reverses the inequality, so
        \begin{equation*}
            (x,r)\longmapsto -\varepsilon^{n/\eta}F^{-1} \bigl(\mu(x+\varepsilon rB^n)\bigr)
        \end{equation*}
    is convex. Letting \(\varepsilon\downarrow0\) and dividing by a positive constant, we obtain joint convexity of
        \begin{equation*}
            P(x,r)=\rho(x)^{-1/\eta}r^{-n/\eta},\qquad r>0.
        \end{equation*}
    Put \(k=n/\eta\). Partial minimization of the jointly convex function \(P(x,r)+kr\) over \(r>0\) preserves convexity. Moreover,
        \begin{equation*}
            \inf_{r>0}\{\rho(x)^{-1/\eta}r^{-k}+kr\} = (k+1)\rho(x)^{-1/(n+\eta)}.
        \end{equation*}
    This proves \eqref{eq:necessary-density-convexity}.
\end{proof}

We apply the convexity condition just proved to profiles with the expansion \eqref{eq:necessary-profile-asymptotic}. A comparison at $R$, $2R$, and $3R$ will give the restriction on its second term without differentiating the remainder.

\begin{proof}[Proof of Theorem~\ref{thm:necessary-profile-tail}]
    Put $q=n+\eta$. Integrating the leading term in \eqref{eq:necessary-profile-asymptotic} over the transverse variables and then over $x_1>t$ gives
        \begin{equation*}
            1-F(t)\sim A t^{-\eta}, \qquad A= \frac{C}{Z_\varphi} \frac{\pi^{(n-1)/2}\Gamma((\eta+1)/2)}{\eta\,\Gamma((n+\eta)/2)} >0.
        \end{equation*}
    Suppose that \eqref{eq:E} holds. By Theorem~\ref{prop:necessary-small-balls}, the function $x\mapsto\rho(x)^{-1/q}$ is convex. Its restriction to a ray, followed by multiplication by a positive constant, shows that
        \begin{equation*}
            G(r):=\varphi(r^2)^{-1/q}
        \end{equation*}
    is convex on $[0,\infty)$. In particular,
        \begin{equation*}
            G(2R)-\frac{G(R)+G(3R)}2\le0 \qquad(R>0).
        \end{equation*}
    If $\theta=1/2$, the claimed condition is automatic. Assume therefore that $\theta\ne1/2$, and put $\beta=1-2\theta$. The profile expansion gives
        \begin{equation*}
            G(r)=C^{-1/q} \left(r-\frac bq r^\beta+o(r^\beta)\right).
        \end{equation*}
    The linear terms cancel in the midpoint difference, so
        \begin{equation*}
            G(2R)-\frac{G(R)+G(3R)}2 = -\frac{C^{-1/q}b}{q}R^\beta \left(2^\beta-\frac{1+3^\beta}{2}\right) +o(R^\beta).
        \end{equation*}
    To determine the sign of the bracket, we consider the auxiliary function $t\mapsto t^\beta$. If $0<\theta<1/2$, then $0<\beta<1$ and this function is strictly concave. Hence
        \begin{equation*}
            2^\beta-\frac{1+3^\beta}{2}>0.
        \end{equation*}
    If $\theta>1/2$, then $\beta<0$ and the auxiliary function is strictly convex. Hence
        \begin{equation*}
            2^\beta-\frac{1+3^\beta}{2}<0.
        \end{equation*}

    Assume, to the contrary, that $b(2\theta-1)>0$. In the first case, this means $b<0$, while in the second case it means $b>0$. Thus the coefficient
        \begin{equation*}
            -\frac{C^{-1/q}b}{q} \left(2^\beta-\frac{1+3^\beta}{2}\right)
        \end{equation*}
    is positive in both cases. Since the remainder is $o(R^\beta)$, for all sufficiently large $R$ we obtain
        \begin{equation*}
            G(2R)-\frac{G(R)+G(3R)}2>0,
        \end{equation*}
    contrary to the convexity of $G$. Therefore $b(2\theta-1)\le0$.
\end{proof}

% \begin{proof}[Proof of Theorem~\ref{thm:necessary-profile-tail}]
%     Integrating the leading term in \eqref{eq:necessary-profile-asymptotic} over transverse directions and then over \(x_1>t\) shows that \(1-F(t)\sim A t^{-\eta}\) for some \(A>0\). Indeed,
%         \begin{equation*}
%             A=\frac{C}{Z_\varphi} \frac{\pi^{(n-1)/2}\Gamma((\eta+1)/2)}{\eta\,\Gamma((n+\eta)/2)}.
%         \end{equation*}
%     Suppose that \eqref{eq:E} holds. Theorem~\ref{prop:necessary-small-balls} then implies convexity of $G(r):=\varphi(r^2)^{-1/(n+\eta)}$ on $[0,\infty)$. Assume, to the contrary, that $b(2\theta-1)>0$. Writing $q=n+\eta$ and expanding \eqref{eq:necessary-profile-asymptotic} gives
%         \begin{equation*}
%             G(r)=C^{-1/q}\left(r-\frac bq r^{1-2\theta} + o(r^{1-2\theta})\right).
%         \end{equation*}
%     Hence, as \(R\to\infty\),
%         \begin{align*}
%             G(2R)-\frac{G(R)+G(3R)}2 ={}&-\frac{C^{-1/q}b}{q}R^{1-2\theta} \left[2^{1-2\theta} -\frac{1+3^{1-2\theta}}2\right] +o(R^{1-2\theta}).
%         \end{align*}
%     The expression in square brackets has the opposite sign to \(2\theta-1\), except at \(2\theta=1\), when it vanishes. Since we assumed $b(2\theta-1)>0$, the leading term is strictly positive for all sufficiently large $R$, contradicting convexity of $G$. Hence $b(2\theta-1)\le0$, as required. %The proof does not differentiate the asymptotic remainder. } \tz{ <-- to zdanie już było, przed zdaniem ''Proof of theorem'' 4 na stronie 18}
% \end{proof}

The midpoint difference in the preceding proof need not have the same sign for small and large $R$. For example, if $\varphi(0)>0$ and the right derivative $\varphi'(0)<0$ exists, then, with $q=n+\eta$,
    \begin{equation*}
        G(2R)-\frac{G(R)+G(3R)}2 = \frac{\varphi'(0)}{q\,\varphi(0)^{1+1/q}}R^2+o(R^2)<0 \qquad(R\downarrow0).
    \end{equation*}
This is the sign required by convexity for this particular midpoint test. The proof above only needs the opposite sign for some large $R$ to rule out global convexity. A correct sign for small $R$ does not by itself prove a local version of \eqref{eq:E}.

We now prove the exact parameter range in Corollary~\ref{cor:necessary-mixture-iff}. Sufficiency was established in Section~\ref{sec:examples}; necessity follows from the tail condition.

\begin{proof}[Proof of Corollary~\ref{cor:necessary-mixture-iff}]
    Sufficiency is Corollary~\ref{cor:rs-adjacent-mixtures}. For necessity, observe that
        \begin{equation*}
            \varphi_\varepsilon(s)=s^{-p} \left(1+(\varepsilon-p\tau)s^{-1}+O(s^{-2})\right).
        \end{equation*}
    If \(\varepsilon>p\tau\), apply Theorem~\ref{thm:necessary-profile-tail} with \(\eta=2p-n>0\), \(\theta=1\), and \(b=\varepsilon-p\tau>0\). Thus the inequality fails. The endpoint \(\varepsilon=p\tau\) is included in the positive part.
\end{proof}

\begin{remark}
    The next section applies Theorem~\ref{thm:necessary-profile-tail} to stable measures. For $0<\alpha<2$, $\alpha\ne1$, their radial density expansion has the form \eqref{eq:necessary-profile-asymptotic} with $\eta=\alpha$, $\theta=\alpha/2$, and $b$ of the same sign as $\alpha-1$. Thus $b(2\theta-1)>0$, which excludes \eqref{eq:E}. We then use the one-dimensional marginal to construct the bounded boxes asserted in Corollary~\ref{thm:stable-classification}.
\end{remark}

\section{Classification of rotationally invariant stable measures}
\label{sec:stable-classification}

We prove Corollary~\ref{thm:stable-classification} in three steps. First, we recall the Cauchy and Gaussian cases. Next, we apply Theorem~\ref{thm:necessary-profile-tail} and the classical density expansions \cite{Kalinauskaite1970,Zaigraev2006} to exclude the other indices. Finally, we construct bounded symmetric boxes for which the inequality fails. This last step uses the one-dimensional expansion in \cite{MatsuiTakemura2006} and proves the additional assertion about the form of the counterexamples.

\begin{proof}[Proof of Corollary~\ref{thm:stable-classification}]
    For $\alpha=1$, the measure in \eqref{eq:stable-characteristic} is the standard isotropic Cauchy law, so Corollary~\ref{cor:cauchy} applies. For $\alpha=2$, its characteristic function is $e^{-|\xi|^2}$, so the measure is $N(0,2I_n)$. Theorem~\ref{thm:gaussian}, with $\beta=1/4$, therefore applies. This proves the sufficient direction.

    \medskip

    We now exclude the remaining indices. Let $0<\alpha<2$ and $\alpha\ne1$. The isotropic stable density has the standard expansion, as $|x|\to\infty$,
        \begin{equation}
            \label{eq:stable-radial-expansion}
            \rho_{n,\alpha}(x) =c_{n,1}(\alpha)|x|^{-n-\alpha} +c_{n,2}(\alpha)|x|^{-n-2\alpha} +o(|x|^{-n-2\alpha}),
        \end{equation}
    where, in the normalization $e^{-|\xi|^\alpha}$,
        \begin{align*}
            c_{n,1}(\alpha) &=-\frac{2^\alpha\Gamma((n+\alpha)/2)}{\pi^{n/2}\Gamma(-\alpha/2)}>0,\\
            c_{n,2}(\alpha) &=\frac{2^{2\alpha-1}\Gamma((n+2\alpha)/2)}{\pi^{n/2}\Gamma(-\alpha)}.
        \end{align*}
    These are the first two nonzero terms of the classical multidimensional stable-density expansion \cite{Kalinauskaite1970,Zaigraev2006}. Their signs can also be read from the inverse Fourier transform of the homogeneous powers $-|\xi|^\alpha$ and $\tfrac12|\xi|^{2\alpha}$. Since $\Gamma(-\alpha)<0$ for $0<\alpha<1$ and $\Gamma(-\alpha)>0$ for $1<\alpha<2$, the coefficient $c_{n,2}(\alpha)$ has the sign of $\alpha-1$. Writing $s=|x|^2$, the squared-radius profile in \eqref{eq:stable-radial-expansion}, up to a positive normalization, therefore satisfies the hypotheses of Theorem~\ref{thm:necessary-profile-tail} with
        \begin{equation*}
            \eta=\alpha,\qquad\theta=\frac\alpha2,\qquad b=\frac{c_{n,2}(\alpha)}{c_{n,1}(\alpha)}.
        \end{equation*}
    Thus $b(2\theta-1)=b(\alpha-1)>0$. The necessary profile-tail condition fails, proving that \eqref{eq:E} cannot hold. Together with the cases $\alpha=1$ and $\alpha=2$, this proves the equivalence in \eqref{eq:stable-iff}. It remains to construct the bounded-box counterexamples.

    \medskip

    Fix $\alpha\in(0,2)\setminus\{1\}$ for the rest of the proof. The first coordinate has characteristic function
        \begin{equation*}
            \widehat\mu_{n,\alpha}(t e_1)=e^{-|t|^\alpha}.
        \end{equation*}
    Its distribution is independent of $n$ and symmetric about the origin. It has a continuous, strictly positive density. Its distribution function $F_\alpha$ is therefore continuous and strictly increasing, and
        \begin{equation*}
            F_\alpha(-t)=1-F_\alpha(t),\qquad t\in\R.
        \end{equation*}
    For $r>0$, put
        \begin{equation*}
            P_r=[-r,r]\times\R^{n-1},\qquad Q_\alpha(r)=F_\alpha^{-1}\!\left(2F_\alpha(r)-1\right).
        \end{equation*}
    Since
        \begin{equation*}
            \mu_{n,\alpha}(P_r)=2F_\alpha(r)-1, \qquad (1-\lambda)P_r+\lambda P_s = P_{(1-\lambda)r+\lambda s},
        \end{equation*}
    property \eqref{eq:E} would imply that $Q_\alpha$ is concave on $(0,\infty)$. In particular, it would imply
        \begin{equation}
            \label{eq:stable-required-midpoint}
            Q_\alpha(2R) -\frac{Q_\alpha(R)+Q_\alpha(3R)}2\ge0 \qquad(R>0).
        \end{equation}
    We will prove that the left-hand side is strictly negative for every sufficiently large $R$.

    \medskip

    To estimate the midpoint difference, let $f_\alpha=F_\alpha'$ and $\overline F_\alpha=1-F_\alpha$. In the normalization $e^{-|t|^\alpha}$, the classical stable-density expansion is
        \begin{align}
            f_\alpha(r) ={}&\frac{\Gamma(1+\alpha)\sin(\pi\alpha/2)}{\pi} r^{-1-\alpha}\notag\\
            &-\frac{\Gamma(1+2\alpha)\sin(\pi\alpha)}{2\pi}r^{-1-2\alpha}+o(r^{-1-2\alpha}),\qquad r\to\infty.
            \label{eq:stable-density-asymptotics}
        \end{align}
    See Matsui and Takemura~\cite[(2.9)]{MatsuiTakemura2006}. For $\alpha<1$, the underlying series converges; for $1<\alpha<2$, it is an asymptotic expansion. Only the two displayed terms and the stated remainder are needed here.

    Integrating \eqref{eq:stable-density-asymptotics} gives
        \begin{equation}
            \label{eq:stable-tail-asymptotics}
            \overline F_\alpha(r) =a_\alpha r^{-\alpha} \bigl(1+b_\alpha r^{-\alpha}+o(r^{-\alpha})\bigr),
        \end{equation}
    where
        \begin{equation*}
            a_\alpha=\frac{\Gamma(\alpha)\sin(\pi\alpha/2)}\pi>0, \qquad b_\alpha=-\frac{\Gamma(2\alpha)}{\Gamma(\alpha)}\cos\frac{\pi\alpha}{2}.
        \end{equation*}
    Indeed, if $h(r)=o(r^{-1-2\alpha})$, then for every $\eta>0$ and all sufficiently large $R$,
        \begin{equation*}
            \left|\int_R^\infty h(r)\dd r\right|\le\eta\int_R^\infty r^{-1-2\alpha}\dd r =\frac{\eta}{2\alpha}R^{-2\alpha}.
        \end{equation*}
    This justifies the integrated remainder in \eqref{eq:stable-tail-asymptotics}. The signs of the second coefficient are
        \begin{equation*}
            b_\alpha<0\quad(0<\alpha<1),\qquad b_\alpha>0\quad(1<\alpha<2).
        \end{equation*}

    \medskip

    To determine the sign of the midpoint difference, we need the second term of $Q_\alpha(r)$, since the leading linear term cancels. By the definition of $Q_\alpha$, we have
        \begin{equation}
            \label{eq:stable-tail-identity}
            \overline F_\alpha(Q_\alpha(r)) = 2\overline F_\alpha(r).
        \end{equation}
    In particular, $Q_\alpha(r)\to\infty$. The leading term in \eqref{eq:stable-tail-asymptotics} gives
        \begin{equation*}
            \frac{Q_\alpha(r)}r\longrightarrow c_\alpha:=2^{-1/\alpha}.
        \end{equation*}
    Consequently,
        \begin{equation*}
            Q_\alpha(r)^{-\alpha} =2r^{-\alpha}+o(r^{-\alpha}).
        \end{equation*}
    Using the full expansion in \eqref{eq:stable-tail-identity} now yields
        \begin{align*}
            \left(\frac{Q_\alpha(r)}{c_\alpha r}\right)^\alpha &=\frac{1+b_\alpha Q_\alpha(r)^{-\alpha} + o(Q_\alpha(r)^{-\alpha})}{1+b_\alpha r^{-\alpha}+o(r^{-\alpha})}\\
            &=1+b_\alpha r^{-\alpha}+o(r^{-\alpha}).
        \end{align*}
    Taking the positive $\alpha$th root, we obtain
        \begin{equation}
            \label{eq:stable-inverse-asymptotics}
            Q_\alpha(r) = c_\alpha r+d_\alpha r^{1-\alpha} +o(r^{1-\alpha}), \qquad d_\alpha=\frac{c_\alpha b_\alpha}{\alpha}.
        \end{equation}
    No derivative of an asymptotic remainder has been taken.

    \medskip
    
    We now compute the midpoint difference. Define
        \begin{equation*}
            \mathcal D_\alpha(R) :=Q_\alpha(2R) -\frac{Q_\alpha(R)+Q_\alpha(3R)}2.
        \end{equation*}
    The linear terms in \eqref{eq:stable-inverse-asymptotics} cancel, leaving
        \begin{equation}
            \label{eq:stable-midpoint-asymptotics}
            \mathcal D_\alpha(R) = d_\alpha \left(2^{1-\alpha}-\frac{1+3^{1-\alpha}}2\right) R^{1-\alpha}+o(R^{1-\alpha}).
        \end{equation}
    If $0<\alpha<1$, then $d_\alpha<0$ and the bracket is positive, since $r\mapsto r^{1-\alpha}$ is strictly concave. If $1<\alpha<2$, then $d_\alpha>0$ and the bracket is negative, since $r\mapsto r^{1-\alpha}$ is strictly convex. Thus the leading coefficient in \eqref{eq:stable-midpoint-asymptotics} is strictly negative in both cases. There exists $R_0=R_0(\alpha)>0$ such that
        \begin{equation*}
            \mathcal D_\alpha(R)<0\qquad(R\ge R_0).
        \end{equation*}
    This contradicts \eqref{eq:stable-required-midpoint} for the corresponding pair of slabs. The inequality in \eqref{eq:stable-required-midpoint} would have to hold for every $R>0$, so its failure for large $R$ is sufficient. We do not claim that $Q_\alpha$ is globally convex; the argument does not exclude concavity on other intervals.

    \medskip
    
    Finally, we truncate these slabs to obtain bounded boxes. Fix any $R\ge R_0$. For $M>0$ and $r>0$, set
        \begin{equation*}
            B_{r,M}=[-r,r]\times[-M,M]^{n-1}.
        \end{equation*}
    These are bounded boxes symmetric with respect to every coordinate hyperplane, and
        \begin{equation*}
            \frac12B_{R,M}+\frac12B_{3R,M}=B_{2R,M}.
        \end{equation*}
    For each fixed $r\in\{R,2R,3R\}$, the boxes increase to $P_r$ as $M\to\infty$, so
        \begin{equation*}
            \mu_{n,\alpha}(B_{r,M}) \longrightarrow 2F_\alpha(r)-1\in(0,1).
        \end{equation*}
    By continuity of $F_\alpha^{-1}$ at the three limiting values in $(0,1)$,
    % By continuity of $F_\alpha^{-1}$ at these three limiting masses,
        \begin{align*}
            &F_\alpha^{-1}\!\left(\mu_{n,\alpha}(B_{2R,M})\right)\\
            &\quad-\frac12F_\alpha^{-1}\!\left(\mu_{n,\alpha}(B_{R,M})\right)-\frac12F_\alpha^{-1}\!\left(\mu_{n,\alpha}(B_{3R,M})\right)\longrightarrow\mathcal D_\alpha(R)<0.
        \end{align*}
    Hence the inequality also fails for all sufficiently large finite $M$. This proves the bounded-box assertion and completes the proof of Corollary~\ref{thm:stable-classification}.
\end{proof}

\begin{remark}
    The index $\alpha$ is fixed throughout the asymptotic argument. We do not claim that the same $R_0$ works when $\alpha$ tends to $1$ or $2$. The result also holds after a common change of scale: for characteristic function $e^{-(\sigma|\xi|)^\alpha}$, $\sigma>0$, the inverse marginal distribution function is $\sigma F_\alpha^{-1}$. The assertion then follows by dilating the sets.
\end{remark}

\section*{Declaration of Generative AI Use}

During the preparation of this manuscript, the authors utilized ChatGPT (OpenAI) as an auxiliary tool for language editing, stylistic refinement, and technical formatting. The model was also used to review the mathematical arguments for correctness and logical consistency. 

All AI-assisted content, mathematical formulations, and arguments were critically reviewed, verified, and edited by the authors to ensure mathematical rigor and accuracy. The authors take full responsibility for the scientific integrity and final content of this work.

\end{document}